\documentclass[a4paper,pdftex,reqno,10pt]{amsart}

\allowdisplaybreaks[1]

\usepackage{geometry}

\usepackage{graphicx}
\usepackage{enumerate}
\usepackage{courier}
\usepackage{times}
\usepackage{amsmath,amssymb}

\usepackage{hyperref}

\theoremstyle{plain}
\newtheorem{Theorem}{Theorem}[section]
\newtheorem{Proposition}[Theorem]{Proposition}

\newtheorem{Corollary}[Theorem]{Corollary}
\newtheorem{Lemma}[Theorem]{Lemma}

\newenvironment{Proof}
{\begin{trivlist}\item[]{{\sc Proof.}}}{\hfill{$\square$}\noindent\end{trivlist}}

\theoremstyle{definition}
\newtheorem{Definition}[Theorem]{Definition}

\newtheorem{Example}[Theorem]{Example}

\theoremstyle{remark}

\newcommand{\cM}{\mathcal{M}}
\newcommand{\cI}{\mathcal{I}}
\newcommand{\cF}{\mathcal{F}}
\newcommand{\cH}{\mathcal{H}}
\newcommand{\cK}{\mathcal{K}}
\newcommand{\cP}{\mathcal{P}}
\newcommand{\cS}{\mathcal{S}}
\newcommand{\PG}{\operatorname{PG}}
\newcommand{\N}{\mathbb{N}}
\newcommand{\Z}{\mathbb{Z}}
\newcommand{\Q}{\mathbb{Q}}

\newcommand{\F}{\mathbb{F}}

\begin{document}


\title{Non-projective two-weight codes}

\date{}

\author{Sascha Kurz}
\address{Sascha Kurz, University of Bayreuth and Friedrich-Alexander-University Erlangen-Nuremberg, Germany}
\email{sascha.kurz@uni-bayreuth.de and sascha.kurz@fau.de}

\begin{abstract}
It has been known since the 1970’s that the difference of the non-zero weights of a projective $\F_q$-linear two-weight has to be a power of the characteristic of the 
underlying field. Here we study non-projective two-weight codes and e.g.\ show the same result under mild extra conditions. For small dimensions we give exhaustive
enumerations of the feasible parameters in the binary case.

\end{abstract}
\maketitle

\section{Introduction}
\noindent
It has been known since the 1970’s that the two non-zero weights of a projective $\F_q$-linear two-weight code $C$ can be written as $w_1=up^t$ and $w_2=(u+1)p^t$, where $u\in\mathbb{N}_{\ge 1}$ and 
$p$ is the characteristic of the underlying finite field $\F_q$, see \cite[Corollary 2]{delsarte1972weights}. So, especially the weight difference $w_2-w_1$ is a power of the characteristic $p$. Here we 
want to consider $\F_q$-linear two-weight codes $C$ with non-zero weights $w_1<w_2$ which are not necessarily projective. In \cite{brouwer2021two} it was observed that if $w_2-w_1$ is not a power of
the characteristic $p$, then the code $C$ has to be non-projective, which settles a question in \cite{luo2017construction}. Here we prove the stronger statement that $C$ is repetitive, i.e., 
$C$ is the $l$-fold repetition of a smaller two-weight code $C'$, where $l$ is the largest factor of $w_2-w_1$ that is coprime to the field size $q$, if $C$ does not have full 
support, c.f.\ \cite{boyvalenkov2021two}. Moreover, if a two-weight code $C$ is non-repetitive and does not have full support, then its two non-zero weights can be written as $w_1=up^t$ and 
$w_2=(u+1)p^t$, where again $p$ is the characteristic of the underlying finite field $\F_q$, see Theorem~\ref{main_thm_3}

Constructions for projective two-weight codes can be found in the classical survey paper \cite{calderbank1986geometry}. Many research papers considered these objects since they e.g.\ yield 
strongly regular graphs and we refer to \cite{brouwer2022strongly} for a corresponding monograph on srgs. For a few more recent papers on constructions for projective two-weight codes we refer  e.g.\ 
to \cite{heng2021family,kohnert2007constructing,pavese2015geometric,zhu2023two}. In e.g.\ \cite{pavese2015geometric} the author uses geometric language and speaks of constructions for two-character sets, 
i.e., sets of points in a projective space $\PG(k-1,q)$ with just two different hyperplane multiplicities, call them $s$ and $t$. In general each (full-length) linear code is in one-to-one correspondence
to a (spanning) multiset of points in some projective space $\PG(k-1,q)$. Here we will also mainly use the geometric language and consider a few general constructions for two-character multisets of points 
corresponding to two-weight codes (possibly non-projective). For each subset $\overline{\cH}$ of hyperplanes in $\PG(k-1,q)$ we construct a multiset of points $\cM(\overline{\cH})$ such 
that all hyperplanes $H\in\overline{\cH}$ have the same multiplicity, say $s$, and also all other hyperplanes $H\notin\overline{\cH}$ have the same multiplicity, say $t$. Actually, we 
characterize the full set of such multisets with at most two different hyperplane multiplicities given $\overline{\cH}$, see Theorem~\ref{main_theorem} and Theorem~\ref{main_thm_2}. Using this correspondence we have 
classified all two-weight codes up to symmetry for small parameters. For projective two-weight codes such enumerations can be found in \cite{bouyukliev2006projective}.  

Brouwer and van Eupen give a correspondence between arbitrary projective codes and arbitrary two-weight codes via the so-called BvE-dual transform. The correspondence can be said to be 1-1, even 
though there are choices for the involved parameters to be made in both directions.  In \cite{brouwer1997correspondence} the dual transform was e.g.\ applied to the unique projective $[16, 5, 9]_3$-code. 
For parameters $\alpha = \tfrac{1}{3}$, $\beta=-3$ the result is a $[69, 5, 45]_3$ two-weight code and for $\alpha=-\tfrac{1}{3}$, $\beta=5$ the results is a $[173, 5, 108]_3$ two-weight code. This resembles 
the fact that we have some freedom when constructing a two-weight code from a given projective code, e.g.\ we can take complements or add simplex codes of the same dimension. Our obtained results 
may be rephrased in the language of the BvE-dual transform by restricting to a canonical choice of the involved parameters. For further literature on the dual transform see e.g.\
\cite{bouyukliev2023dual,bouyukliev2009classification,brouwer1997correspondence,takenaka2008optimal}. For a variant that is rather close to our presentation we refer to \cite{bouyuklieva2017dual}.

With respect to further related literature we remark that a special subclass of (non-projective) two-weight codes was completely characterized in \cite{jungnickel2018classification}. 
A conjecture by Vega \cite{vega2012weight} states that all two-weight cyclic codes are the {\lq\lq}known{\rq\rq} ones, c.f.\ \cite{duc2019non}. Another stream of literature considers the problem 
whether all projective two-weight codes that have the parameters of partial $k$-spreads indeed have to be partial $k$-spreads. Those results can be found in papers considering extendability results 
for partial $k$-spreads or classifying minihypers, see e.g.\ \cite{govaerts2003particular}. Several non-projective two-weight codes appear also as minimum length examples for divible minimal 
codes \cite{kurz2023divisible}.\footnote{Minimal codes are linear codes where all non-zero codewords are minimal, i.e., whose support is not properly contained in the support of another codeword.} 
Two-weight codes have also been considered over rings instead of finite fields, see e.g.\ \cite{byrne2008ring}.

The remaining part of this paper is structured as follows. In Section~\ref{sec_preliminaries} we introduce the necessary preliminaries for linear two-weight codes and their geometric counterpart called
two-character multisets in projective spaces. In general multisets of points, corresponding to general linear codes, can be described via so-called characteristic functions and we collect 
some of their properties in Section~\ref{sec_characteristic_functions}. Examples and constructions for two-character multisets are given in Section~\ref{sec_constructions}. In Section~\ref{sec_feasibility_sets} 
we present our main results. We close with enumeration results for two-character multisets in $\PG(k-1,q)$ for small parameters in Section~\ref{sec_enumeration}. We will mainly use geometric 
language and arguments. For the for the ease of the reader we only use elementary arguments and give (almost) all proof details. 
In an appendix we list all candidates for sets of parameters of two-character sets in $\PG(k-1,2)$ for $4\le k\le 12$ and state corresponding constructions from the literature. A few cases are 
left open for $k\in\{11,12\}$.

\section{Preliminaries}
\label{sec_preliminaries}
\noindent
An $[n,k]_q$-code $C$ is a $k$-dimensional subspace of $\F_q^n$, i.e., $C$ is assumed to be $\F_q$-linear. Here $n$ is called the length and $k$ is called the dimension of $C$. Elements $c\in C$ are 
called codewords and the weight $\operatorname{wt}(c)$ of a codeword is given by the number of non-zero coordinates. Clearly, the all-zero vector $\mathbf{0}$ has weight zero and all other codewords 
have a positive integer weight. A two-weight code is a linear code with exactly two non-zero weights. A generator matrix for $C$ is an $k\times n$ matrix $G$ such that its rows span $C$. 
We say that $C$ is of full length if for each index $1\le i\le n$ there exists a codeword $c\in C$ whose $i$th coordinate $c_i$ is non-zero, i.e., all columns of a generator matrix of $C$ are 
non-zero. The dual code $C^\perp$ of $C$ is the $(n-k$)-dimensional code consisting of the vectors orthogonal to all codewords of $C$ w.r.t.\ the inner product $\langle u, v\rangle = \sum_{i=1}^n u_iv_i$. 

Now let $C$ be a full-length $[n,k]_q$-code with generator matrix $G$. Each column $g$ of $G$ is an element of $\F_q^k$ and since $g\neq\mathbf{0}$ we can consider $\langle g\rangle$ as point in the projective 
space $\PG(k-1,q)$. Using the geometric language we call $1$-, $2$-, $3$-, and $(k-1)$-dimensional subspaces of $\F_q^k$ points, lines, planes, and hyperplanes in $\PG(k-1,q)$. Instead of an 
$l$-dimensional space we also speak of an $l$-space. By $\cP$ we denote the 
set of points and by $\cH$ we denote the set of hyperplanes of $\PG(k-1,q)$ whenever $k$ and $q$ are clear from the context. A multiset of points in $\PG(k-1,q)$ is a mapping $\cM\colon\cP\to\N$, i.e., 
to each point $P\in\cP$ we assign its multiplicity $\cM(P)\in\N$. By $\#\cM=\sum_{P\in\cP} \cM(P)$ we denote the cardinality of $\cM$. The support $\operatorname{supp}(\cM)$ is the set 
of all points with non-zero multiplicity. We say that $\cM$ is spanning if the set of points in the support of $\cM$ span $\PG(k-1,q)$. Clearly permuting columns of a generator matrix $G$ or multiplying 
some columns with non-zero elements in $\F_q^\star:=\F_q\backslash\{0\}$ yields an equivalent code. Besides that we get an one-to-one correspondence between full length $[n,k]_q$-codes and spanning
multisets of points $\cM$ in $\PG(k-1,q)$ with cardinality $\#\cM=n$. Moreover, two linear $[n,k]_q$-codes $C$ and $C'$ are equivalent iff their corresponding multisets of points $\cM$ and $\cM'$ are. For 
details we refer e.g.\ to \cite{dodunekov1998codes}. A linear code $C$ is projective iff its corresponding multiset of points satisfies $\cM(P)\in\{0,1\}$ for all $P\in\cP$. We also speak of a set 
of points in this case. The multisets of points with $\cM(P)=0$ for all $P\in\cP$ are called trivial. 

Geometrically, for a non-zero codeword $c\in C$ the set $c\cdot \F_q^\star$ corresponds to a hyperplane $H\in\cH$ and $\operatorname{wt}(c)=\#\cM-\cM(H)$, where we extend the function $\cM$ additively.
i.e., $\cM(S):=\sum_{P\in S} \cM(P)$ for every subset $S\subseteq\cP$ of points. We call $\cM(H)$ the multiplicity of hyperplane $H\in\cH$ and have $\cM(V)=\#\cM$ for the entire ambient space $V:=\cP$.
The number of hyperplanes $\#\cH$, as well as the number of points $\#\cP$, in $\PG(k-1,q)$ is given by $[k]_q:=\tfrac{q^k-1}{q-1}$. A two-character multiset is a multiset of points $\cM$ such that exactly
two different hyperplane multiplicities $\cM(H)$ occur. I.e., a multiset of points $\cM$ is a two-character multiset iff its corresponding code $C$ is a two-weight code. If $\cM$ actually is a set of
points, i.e.\ if we have $\cM(P)\in\{0,1\}$ for all points $P\in\cP$, then we speak of a two-character set. We say that an $[n,k]_q$-code $C$ 
is $\Delta$-divisible if the weights of all codewords are divisible by $\Delta$. A multiset of points $\cM$ is called $\Delta$-divisible if the corresponding linear code is. More directly, a multiset of 
points $\cM$ is $\Delta$-divisible if we have $\cM(H)\equiv\#\cM \pmod\Delta$ for all $H\in\cH$.    

A one-weight code is an $[n,k]_q$-code $C$ such that all non-zero codewords have the same weight $w$. One-weight codes have been completely classified in \cite{bonisoli1984every} and are given by 
repetitions of so-called simplex codes. Geometrically, the multiset of points $\cM$ in $\PG(k-1,q)$ corresponding to a one-weight code $C$ satisfies $\cM(P)=l$ for all $P\in\cP$. I.e., we 
have $\#\cM=n=[k]_q\cdot l$, $\cM(H)=[k-1]_q\cdot l$ for all $H\in \cH$, and $w=\#\cM-\cM(H)=q^{k-1}\cdot l$. We say that a linear $[n,k]_q$-code $C$ is repetitive if it is the $l$-fold repetition of 
an $[n/l,k]_q$-code $C'$, where $l>1$, and non-repetitive otherwise. A given multiset of points $\cM$ is called repeated if its corresponding code is. More directly, a non-trivial multiset of points 
$\cM$ is repeated iff the greatest common divisor of all point multiplicities is larger than one. We say that a multiset of points $\cM$ or its corresponding linear code $C$ has full support 
iff $\operatorname{supp}(\cM)=\cP$, i.e., if $\cM(P)>0$ for all $P\in\cP$. So, for each non-repetitive one-weight code $C$ with length $n$, dimension $k$, and non-zero weight $w$ we have $n=[k]_q$ and $w=q^{k-1}$. 
Each non-trivial one-weight code, i.e., one with dimension at least $1$, has full support. The aim of this paper is to characterize the possible parameters of non-repetitive two-weight codes (with or without 
full support). For the correspondence between $[n,k]_q$-codes and multisets of points $\cM$ in $\PG(k-1,q)$ we have assumed that $\cM$ is spanning. If $\cM$ is not spanning, then there exists a hyperplane
containing the entire support, so that $\cM$ is two-character multiset iff $\cM$ induces a one-character multiset in the span of $\operatorname{supp}(\cM)$ cf.~Proposition~\ref{prop_subspace}. The 
structure of the set of all two-character multisets where the larger hyperplane multiplicity is attained for a prescribed subset of the hyperplanes is considered in Section~\ref{sec_feasibility_sets}.    

\section{Characteristic functions}
\label{sec_characteristic_functions}
\noindent
Fixing the field size $q$ and the dimension $k$ of the ambient space, a multiset of points in $\PG(k-1,q)$ is a mapping $\cM\colon \cP\to\N$. By $\cF$ we denote the $\Q$-vector space 
consisting of all functions $F\colon\cP\to\Q$, where addition and scalar multiplication is defined pointwise. I.e., $(F_1+F_2)(P):=F_1(P)+F_2(P)$ and $(x\cdot F_1)(P):=x\cdot F_1(P)$ for all 
$P\in\cP$, where $F_1,F_2\in\cF$ and $x\in\Q$ are arbitrary. For each non-empty subset $S\subseteq\cP$ the characteristic function $\chi_S$ is defined by $\chi_S(P)=1$ if $P\in S$ and 
$\chi_S(P)$ otherwise. Clearly the set of functions $\chi_P$ for all $P\in\cP$ forms a basis of $\cF$ for ambient space $\PG(k-1,q)$ for all $k\ge 1$. Note that there are no hyperplanes 
if $k=1$ and hyperplanes coincide with points for $k=2$. We also extend the functions $F\in\cF$ additively, i.e., we set $F(S)=\sum_{P\in S} F(P)$ for all $S\subseteq\cP$. Our 
next aim is to show the well-known fact that also the set of functions $\chi_H$ for all hyperplanes $H\in\cH$ forms a basis of $\cF$ for ambient space $\PG(k-1,q)$ for all $k\ge 2$. In other words, 
also $\cM(P)$ can be reconstructed from the $\cM(H)$:
\begin{Lemma}
  \label{lemma_reconstruct_points_from_hyperplanes}
  Let $\cM\in\cF$ for ambient space $\PG(k-1,q)$, where $k\ge 2$. Then, we have
  \begin{equation}
    \cM(P)=\sum_{H\in \cH\,:\, P\in H} \frac{1}{[k-1]_q}\cdot \cM(H)+ \sum_{H\in \cH\,:\, P\notin H} \frac{1}{q^{k-1}}\cdot\left(\frac{1}{[k-1]_q}-1\right)\cdot \cM(H)
  \end{equation}
  for all points $P\in\cP$.
\end{Lemma}
\begin{Proof}
  W.l.o.g.\ we assume $k\ge 3$. 
  Since each point $P'\in\cP$ is contained in $[k-1]_q$ of the $\#\cH=[k]_q$ hyperplanes and each point $P'\neq P$ is contained in $[k-2]_q$ of 
  the $[k-1]_q$ hyperplanes that contain $P$, we have
  $$
     \sum_{H\in\cH\,:\,P\in H} \cM(H) =[k-2]_q\cdot |\cM|+\left([k-1]_q-[k-2]_q\right)\cM(P)= [k-2]_q\cdot |\cM| +q^{k-2}\cM(P)
  $$
  so that
  $$
    \sum_{H\in\cH\,:\,P\in H} \cM(H)-\frac{[k-2]_q}{[k-1]_q}\cdot\sum_{H\in\cH} \cM(H)=q^{k-2} \cM(P)   
  $$  
  using $[k-1]_q\cdot\#\cM=\sum_{H\in\cH} \cM(H)$.
  Thus, we can conclude the stated formula using 
  $$
    \frac{1}{q^{k-2}}\cdot \left(1-\frac{[k-2]_q}{[k-1]_q}\right)=
    \frac{1}{q^{k-2}}\cdot \frac{[k-1]_q-[k-2]_q}{[k-1]_q}=\frac{1}{[k-1]_q}
  $$  
  and
  $$
    -\frac{[k-2]_q}{[k-1]_q \cdot q^{k-2}}=\frac{1-[k-1]_q}{[k-1]_q\cdot q^{k-1}}=\frac{1}{q^{k-1}}\cdot\left(\frac{1}{[k-1]_q}-1\right). 
  $$
\end{Proof}

As an example we state that in $\PG(3-1,2)$ we have 
$$
  \cM(P)=\frac{1}{3} \cdot \sum_{H\in \cH\,:\, P\in H}  \cM(H) -\frac{1}{6} \cdot \sum_{H\in \cH\,:\, P\notin H} \cM(H).
$$

\begin{Lemma}
  \label{lemma_basis_hyperplanes}
  Let $\cM\in\cF$ for ambient space $\PG(k-1,q)$, where $k\ge 2$. Then there exist $\alpha_H\in\mathbb{Q}$ for all hyperplanes $H\in\cH$ such that
  \begin{equation}
    \cM=\sum_{H\in\cH} \alpha_H\cdot\chi_H.
  \end{equation}
  Moreover, the coefficients $\alpha_H$ are uniquely determined by $\cM$.
\end{Lemma}
\begin{Proof}
  From 
  $$
    \sum_{H\in\cH\,:\,P\in H} \chi_H-\frac{[k-2]_q}{[k-1]_q}\cdot\sum_{H\in\cH} \chi_H=q^{k-2} \cdot\chi_P   
  $$
  for each point $P\in\cP$ and 
  $$
    \cM=\sum_{P\in\cP} \cM(P)\cdot\chi_P
  $$
  we conclude the existence of the $\alpha_H\in\mathbb{Q}$. Since the functions $\left(\chi_P\right)_{P\in\cP}$ form a basis of the $\Q$-vector space $\cF$ , 
  which is also generated by the functions $\left(\chi_H\right)_{H\in\cH}$, counting $\#\cP=[k]_q=\#\cH$ yields that also $\left(\chi_H\right)_{H\in\cH}$ forms a basis 
  and the coefficients $\alpha_H$ are uniquely determined by $\cM$. 
\end{Proof}

If $\cM\in\cF$ is given by the representation
$$
  \cM=\sum_{P\in\cP} \alpha_P\cdot\chi_P
$$
with $\alpha_P\in\Q$ we can easily decide whether $\cM$ is a multiset of points. The necessary and sufficient conditions are given by $\alpha_P\in\N$ for all $P\in \cP$ (including the
case of a trivial multiset of points). If a multiset of points is characterized by coefficients $\alpha_H$ for all hyperplanes $H\in\cH$ as in Lemma~\ref{lemma_basis_hyperplanes} then 
some $\alpha_H$ may be fractional or negative. For two-character multisets we will construct a different unique representation involving the characteristic functions $\chi_H$ of
hyperplanes, see Theorem~\ref{main_theorem}.

Let us state a few observations about operations for multisets of points that yield multisets of points again.
\begin{Lemma}
  For two multisets of points $\cM_1,\cM_2$ of $\PG(k-1,q)$ and each non-negative integer $n\in\N$ the functions $\cM_1+\cM_2$ and $n\cdot \cM_1$ are multisets of points of $\PG(k-1,q)$. 
\end{Lemma}
In order to say something about the subtraction of multisets of points we denote the minimum point multiplicity of a multiset of points $\cM$ by $\mu(\cM)$ and the maximum point multiplicity 
by $\gamma(\cM)$. Whenever $\cM$ is clear from the context we also just write $\mu$ and $\gamma$ instead of $\mu(\cM)$ and $\mu(\gamma)$.  
\begin{Lemma}
  Let $\cM_1$ and $\cM_2$ be two multisets of points of $\PG(k-1,q)$. If $\mu(\cM_1)\ge\gamma(\cM_2)$, then $\cM_1-\cM_2$ is a multiset of points of $\PG(k-1,q)$. 
\end{Lemma}
\begin{Definition}
  Let $\cM$ be a multiset of points in $\PG(k-1,q)$. If $l$ is an integer with $l\ge \gamma(\cM)$, then the $l$-complement $\cM^{l-C}$ of $\cM$ is defined by 
  $\cM^{l-C}(P):=l-\cM(P)$ for all points $P\in\cP$.
\end{Definition}
One can easily check that $\cM^{l-C}$ is a multiset of points with cardinality $l\cdot[k]_q-\#\cM$, maximum point multiplicity $\gamma\!\left(\cM^{l-C}\right)=l-\mu(\cM)$, and 
minimum point multiplicity $\mu\!\left(\cM^{l-C}\right)=l-\gamma(\cM)$. Using characteristic functions we can write $\cM^{l-C}=l\cdot\chi_V-\cM$, where $V=\cP$ denotes the set 
of all points of the ambient space.

Given an arbitrary function $\cM\in\cF$ there always exist $\alpha\in\Q\backslash\{0\}$ and $\beta\in\Z$ such that $\alpha\cdot\cM+\beta\cdot\chi_V$ is a multiset of points.

\section{Examples and constructions for two-character multisets}
\label{sec_constructions}
The aim of this section is to list a few easy constructions for two-character multisets of points $\cM$ in $\PG(k-1,q)$. We will always abbreviate $n=\#\cM$ and denote the
two occurring hyperplane multiplicities by $s$ and $t$, where we assume $s>t$ by convention.

\begin{Proposition}
  \label{prop_subspace}
  For integers $1\le l<k$ let $L$ be an arbitrary $l$-space in $\PG(k-1,q)$. Then $\chi_L$ is a two-character set with 
  $n=[l]_q$, $\gamma=1$, $\mu=0$, $s=[l]_q$, and $t=[l-1]_q$. 
\end{Proposition}
Note that for the case $l=k$ we have the one-character set $\chi_V$, which can be combined with any two-character multiset.
\begin{Lemma}
  \label{lemma_addition_and_minus}
  Let $\cM$ be a two-character multiset of points in $\PG(k-1,q)$. Then, for each integer $0\le a\le \mu(\cM)$, each $b\in N$, and each integer $c\ge \gamma(\cM)$ the functions 
  $\cM-a\cdot\chi_V$, $\cM+b\cdot\chi_V$, $b\cdot \cM$, and $c\cdot\chi_V-\cM$ are two-character multisets of points. 
\end{Lemma}

For the first and the fourth construction we also spell out the implications for the parameters of a given two-character multiset:
\begin{Lemma}
  \label{lemma_minus_ambient_space}
  Let $\cM$ be a multiset of points in $\PG(k-1,q)$ such that $\cM(H)\in\{s,t\}$ for every hyperplane $H\in\cH$. If $\cM(P)\ge l$ for every point $P\in\cP$, i.e., $l\le \mu(\cM)$, then 
  $\cM':=\cM-l\cdot \chi_V$ is a multiset of points in $\PG(k-1,q)$ such that $\cM'(H)\in\{s-[k-1]_q\cdot l,t-[k-1]_q\cdot l\}$ for every hyperplane $H\in\cH$.
\end{Lemma}

\begin{Lemma}
  \label{lemma_complement}
  Let $\cM$ be a multiset of points in $\PG(k-1,q)$ such that $\cM(H)\in\{s,t\}$ for every hyperplane $H\in\cH$. If $\cM(P)\le u$, i.e.\ $\le \gamma(\cM)$ for every point $P\in\cP$, then 
  the $u$-complement $\cM':=u\cdot\chi_V-\cM$ of $\cM$ is a multiset of points in $\PG(k-1,q)$ such that $\cM'(H)\in\left\{u[k-1]-s,u[k-1]-t\right\}$ for every hyperplane $H\in\cH$.
\end{Lemma}
 
We can also use two (almost) arbitrary subspaces to construct two-character multisets:
\begin{Proposition}
  \label{prop_two_subspaces}
  Let $a\ge b\ge 1$ and $0\le i\le b-1$ be arbitrary integers, $A$ be an $a$-space and $B$ be an $b$-space with $\dim(A\cap B)=i$ in $\PG(k-1,q)$, where $k=a+b-i$, Then, $\cM=\chi_A+q^{a-b}\cdot \chi_B$ 
  satisfies $\cM(H)\in\{s,t\}$ for all $H\in\cH$, where $s=[a-1]_q+q^{a-b}\cdot[b-1]_q$ and $t=s+q^{a-1}$. If $i=0$, then $\gamma=q^{a-b}$ and $\gamma=q^{a-b}+1$ otherwise. In general, we have 
  $n=[a]_q+q^{a-b}\cdot[b]_q$ and $\mu=0$. 
\end{Proposition}
\begin{Proof}
  For each $H\in\cH$ we have $\cM(H\cap A)\in\left\{[a-1]_q,[a]_q\right\}$ and $\cM(H\cap B)\in\left\{[b-1]_q,[b]_q\right\}$. Noting that we cannot have both $\cM(H\cap A)=[a]_q$ and 
  $\cM(H\cap B)=[b]_q$ we conclude $\cM(H)\in\left\{[a-1]_q+q^{a-b}\cdot [b-1]_q,[a]_q+q^{a-b}\cdot [b-1]_q,[a-1]_q+q^{a-b}\cdot [b]_q\right\}=\{s,t\}$.   
\end{Proof}   

A partial $k$-spread is a set of $k$-spaces in $\PG(v-1,q)$ with pairwise trivial intersection.
\begin{Proposition}
  \label{prop_partial_parallelism}
  Let $\cS_1,\dots,\cS_r$ be a partial parallelism of $\PG(2k-1,q)$, i.e., the $\cS_i$ are partial $k$-spreads that are pairwise disjoint.
  Then
  $$
    \cM=\sum_{i=1}^r\sum_{S\in\cS_i} \chi_S
  $$
  is a two-character multiset of $\PG(2k-1,q)$ with $n=r\cdot[k]_q$ and hyperplane multiplicities $s=r\cdot[k-1]_q$, $t=r\cdot[k-1]_q+q^{k-1}$, where
  $r=\sum_{i=1}^r \left|\cS_i\right|$.
\end{Proposition}  
C.f.\ Example SU2 in \cite{calderbank1986geometry}. Field changes work similarly as explained in \cite[Section 6]{calderbank1986geometry} for two-character sets. 

Based on hyperplanes we can construct large families of two-character multisets:
\begin{Lemma}
  \label{lemma_sums_of_hyperplanes}
  Let $\emptyset\neq\cH'\subsetneq \cH$ be a subset of the hyperplanes of $\PG(k-1,q)$, where $k\ge 3$, then
  \begin{equation}
    \cM=\sum_{H\in\cH'} \chi_H
  \end{equation}
  is a two-character multiset with $n=r[k-1]_q$, $s=r[k-2]_q+q^{k-2}$, and $t=r[k-2]_q$, where $r=\# \cH'$.
\end{Lemma}
By allowing $\cH'$ to be a multiset of hyperplanes, we end up with $(\tau+1)$-character sets, where $\tau$ is the maximum number of occurrences of a hyperplane in $\cH'$.

Applying Lemma~\ref{lemma_minus_ambient_space} yields:
\begin{Lemma}
  Let $\emptyset\neq\cH'\subsetneq \cH$ be a subset of the hyperplanes of $\PG(k-1,q)$, where $k\ge 3$. If each point $P\in\cP$ is contained in at 
  least $\mu\in N$ elements of $\cH'$, then
  \begin{equation}
    \cM=\sum_{H\in\cH'} \chi_H\,-\,\mu\cdot \chi_V
  \end{equation}
  is a two-character multiset with $n=r[k-1]_q-\mu[k]_q$, $s=r[k-2]_q+q^{k-2}-\mu[k-1]_q$ and $t=r[k-2]_q-\mu[k-1]_q$, where $r=|\cH'|$.
\end{Lemma}

In some cases we obtain two-character multisets where all point multiplicities have a common factor $g>1$. Here we can apply the following general construction:
\begin{Lemma}
  \label{lemma_gcd}
  Let $\cM$ be a multiset of points in $\PG(k-1,q)$ such that $\cM(H)\in\{s,t\}$ for every hyperplane $H\in\cH$. If $\cM(P)\equiv 0\pmod g$ for every point $P\in\cP$, then 
  $\cM':=\tfrac{1}{g}\cdot \cM$ is a multiset of points in $\PG(k-1,q)$ such that $\cM'(H)\in\left\{\tfrac{1}{g}\cdot s,\tfrac{1}{g}\cdot t\right\}$ for every hyperplane $H\in\cH$. 
  Moreover, we have $\#\cM'=\tfrac{1}{g}\cdot\#\cM$, $\mu(\cM')=\tfrac{1}{g}\cdot\mu(\cM)$, and $\gamma(\cM')=\tfrac{1}{g}\cdot\gamma(\cM)$. 
\end{Lemma}

Interestingly enough, it will turn out that we can construct all two-character multisets by combining Lemma~\ref{lemma_sums_of_hyperplanes} with Lemma~\ref{lemma_addition_and_minus} and 
Lemma~\ref{lemma_gcd}, see Theorem~\ref{main_theorem} and Theorem~\ref{main_thm_2}. 

\section{Geometric duals and sets of feasible parameters for two-character multisets}
\label{sec_feasibility_sets}
\noindent
To each two-character multiset $\cM$ in $\PG(k-1,q)$, i.e., $\left\{\cM(H)\,:\, H\in\cH\right\}=\{s,t\}$ for some $s,t\in\N$ we can assign a set of points
$\overline{\cM}$ by using the geometric dual, i.e., interchanging hyperplanes and points. More precisely, fix a non-degenerated billinear form $\perp$ and consider pairs of points and hyperplanes $(P,H)$ that are 
perpendicular w.r.t. to $\perp$. We write $H=P^\perp$ for the geometric dual of a point. We define $\overline{\cM}$ via $\overline{\cM}(P)=1$ iff $\cM(H)=s$, where $H=P^\perp$, and $\overline{\cM}(P)=0$ otherwise, 
i.e., if $\cM(H)=t$.\footnote{A generalization of the notion of the geometric dual has been introduced by Brouwer and van Eupen \cite{brouwer1997correspondence} for linear codes and formulated for 
multisets of points by Dodunekov and Simonis \cite{dodunekov1998codes}.} Of course we have some freedom how we order $s$ and $t$. So, we may also write $\overline{\cM}(P)=\left(\cM(H)-t\right)/(s-t)\in\{0,1\}$ for all $P\in \cP$, where $H=P^\perp$. Noting 
the asymmetry in $s$ and $t$ we may also interchange the role of $s$ and $t$ or replace $\overline{\cM}$ by its complement. Note that in principle several multisets of points with two hyperplane 
multiplicities can have the same corresponding point set $\overline{\cM}$.

For the other direction we can start with an arbitrary set of points $\overline{\cM}$, i.e., $\overline{\cM}(P)\in\{0,1\}$ for all $P\in\cP$. The multiset of points with two hyperplane multiplicities 
$\cM$ is then defined via $\cM(H)=s$ if $\overline{\cM}(P)=1$, where $H=P^\perp$, and $\cM(H)=t$ if $\overline{\cM}(P)=0$. I.e., we may set
\begin{equation}  
   \cM(H)=t+(s-t)\cdot\overline{\cM}(H^\perp).
\end{equation} 
While we have $\cM(H)\in\N$ for all $s,t\in \N$, the point multiplicities $\cM(P)$ induced by the hyperplane multiplicities $\cM(H)$, see Lemma~\ref{lemma_reconstruct_points_from_hyperplanes}, are not 
integral or non-negative in general. For suitable choices of $s$ and $t$ they are, for others they are not. 

\begin{Definition}
  Let $\overline{\cM}$ be a set of points in $\PG(k-1,q)$. By $\mathbb{L}(\overline{\cM})\subseteq\N^2$ we denote the set of all pairs $(s,t)\in\N^2$ with $s\ge t$ such that a multiset of points $\cM$ in 
  $\PG(k-1,q)$ exists with $\cM(H)=s$ if $\overline{\cM}(H^\perp)=1$ and $\cM(H)=t$ if $\overline{\cM}(H^\perp)=0$ for all hyperplanes $H\in\cH$. 
\end{Definition}

Directly from Lemma~\ref{lemma_addition_and_minus} we can conclude:
\begin{Lemma}
  Let $\overline{\cM}$ be a set of points in $\PG(k-1,q)$. If $(s,t)\in \mathbb{L}(\overline{\cM})$, then we have
  \begin{equation}
    \langle (s,t)\rangle_{\N}+\left\langle \left([k-1]_q,[k-1]_q\right)\right\rangle_{\N}=\left\{\left(us+v[k-1]_q,ut+v[k-1]_q\right)\,:\,u,v\in\N\right\}\subseteq \mathbb{L}(\overline{\cM}). 
  \end{equation}
\end{Lemma}  

Before we study the general structure of $\mathbb{L}(\overline{\cM})$ and show that it can generated by a single element $(s_0,t_0)$ in the above sense, we consider all non-isomorphic examples
in $\PG(3-1,2)$ (ignoring the constraint $s\ge t$). 

\begin{Example}
  \label{ex_one_line}
  Let $\cM$ be a multiset of points in $\PG(2,2)$ uniquely characterized by $\cM(L)=s\in\N$ for some line $L$ and $\cM(L')=t\in\N$ for all other lines $L'\neq L$. 
  For each point $P\in L$ we have 
  \begin{equation}
    \cM(P) = \frac{s+2t}{3}-\frac{4t}{6}=\frac{s}{3}
  \end{equation} 
  and for each point $Q\notin L$ we have
  \begin{equation}
    \cM(Q) = \frac{3t}{3}-\frac{s+3t}{6}=\frac{3t-s}{6}.
  \end{equation} 
  Since $\cM(P),\cM(Q)\in\N$ we set $x:=\cM(P)=\tfrac{s}{3}$ and $y:=\cM(Q)=\tfrac{3t-s}{6}$, so that $s=3x$ and $t=2y+x$. With this we have $n=3x+4y$, $\gamma=\max\{x,y\}$, and 
  $s-t=2(x-y)$. If $x\ge y$, then we can write $\cM=y\cdot \chi_E+(x-y)\cdot\chi_L$. If $x\le y$, then we can write $\cM=y\cdot \chi_E-(y-x)\cdot\chi_L$. 
\end{Example}

For Example~\ref{ex_one_line} the set of all feasible $(s,t)$-pairs assuming $s\ge t$ is given by $\left\langle(3,1)\right\rangle_{\N}+\left\langle(3,3)\right\rangle_{\N}$. If we assume 
$t\ge s$, then the set of feasible $(s,t)$-pairs is given by $\left\langle(0,2)\right\rangle_{\N}+\left\langle(3,3)\right\rangle_{\N}$. The vector $(0,2)$ can be computed from $(3,1)$ 
by computing a suitable complement. 

Due to Lemma~\ref{lemma_minus_ambient_space} we can always assume the existence of a point of multiplicity $0$ as a normalization. So, in Example~\ref{ex_one_line} we may assume $x=0$ or $y=0$, 
so that $\cM=y\cdot \chi_E-y\cdot \chi_L$ or $\cM=x\cdot\chi_L$.

Due to Lemma~\ref{lemma_gcd} we can always assume that the greatest common divisor of all point multiplicities is $1$ as a normalization (excluding the degenerated case of an empty multiset of points). 
Applying both normalizations to the multisets of points in Example~\ref{ex_one_line} leaves the two possibilities $\chi_L$ and $\chi_E-\chi_L$, i.e., point sets. 
 
Due to Lemma~\ref{lemma_complement} we always can assume $\#\cM \le \gamma(\cM)\cdot[k]_q/2$. Applying also the third normalization to the multisets of points in Example~\ref{ex_one_line} leaves only the 
possibility $\chi_L$, i.e., a subspace construction, see Proposition~\ref{prop_subspace}, where $s=3$, $t=1$, $n=3$, and $s-t=2$.

\begin{Example}
  \label{ex_two_lines}
  Let $\cM$ be a multiset of points in $\PG(2,2)$ uniquely characterized by $\cM(L_1)=\cM(L_2)=s\in\N$ for two different lines $L_1,L_2$ and $\cM(L')=t\in\N$ for all other lines 
  $L'\notin\left\{L_1,L_2\right\}$. For $P:=L_1\cap L_2$ we have 
  \begin{equation}
    \cM(P) = \frac{2s+t}{3}-\frac{4t}{6}=\frac{2s-t}{3},
  \end{equation} 
  for each point $Q\in \left(L_1\cup L_2\right)\backslash \{P\}$ we have
  \begin{equation}
    \cM(Q) = \frac{s+2t}{3}-\frac{s+3t}{6}=\frac{s+t}{6},
  \end{equation} 
  and for each point $R\notin L_1\cup L_2$ we have
  \begin{equation}
    \cM(R) = \frac{3t}{3}-\frac{2s+2t}{6}=\frac{2t-s}{3}.
  \end{equation} 
  Since $\cM(Q),\cM(R)\in\N$ we set $x:=\cM(Q)=\tfrac{s+t}{6}$ and $y:=\cM(R)=\tfrac{2t-s}{3}$, so that $s=4x-y$ and $t=2x+y$. With this we have $n=6x+7y$ and 
  $s-t=2(x-y)$. Of course we need to have $y\le 2x$ so that $\cM(P)\ge 0$, which implies $s\ge 0$.
  \begin{itemize} 
    \item $\cM(P)=0$: $y=2x$, so that $\cM(P)=0$, $\cM(Q)=x$, $\cM(R)=2x$, and the greatest common divisor of $\cM(P)$, $\cM(Q)$, and $\cM(R)$ is equal to $x$. Thus, we can 
                      assume $x=1$, $y=2$, so that $s=2$, $t=4$, $n=8$, $\gamma=2$, $t-s=2$, and 
                      $\cM=2\chi_E-\chi_{L_1}-\chi_{L_2}$ for two different lines $L_1,L_2$. 
    \item $\cM(Q)=0$: $x=0$, so that also $y=0$ and $\cM$ is the empty multiset of points.
    \item $\cM(R)=0$: $y=0$, $\cM(P)=2x$, $\cM(Q)=x$, so that $\gcd(\cM(P),\cM(Q),\cM(R))=x$ and we can assume $x=1$. With this we have $s=4$, $t=2$, $n=6$, $\gamma=2$, $s-t=2$, and 
                      $\cM=\chi_{L_1}+\chi_{L_2}$ for two different lines $L_1,L_2$.
  \end{itemize}
\end{Example}
So, Example~\ref{ex_two_lines} can be explained by the construction in Proposition~\ref{prop_two_subspaces}.

\begin{Example}
  \label{ex_three_lines_with_common_intersection}
  Let $\cM$ be a multiset of points in $\PG(2,2)$ uniquely characterized by $\cM(L_1)=\cM(L_2)=\cM(L_3)=s\in\N$ for three different lines $L_1,L_2,L_3$ with a common 
  intersection point $P=L_1\cap L_2\cap L_3$ and $\cM(L')=t\in\N$ for all other lines. We have
  \begin{equation}
    \cM(P) = \frac{3s}{3}-\frac{4t}{6}=s-\frac{2t}{3}
  \end{equation} 
  and 
  \begin{equation}
    \cM(Q) = \frac{s+2t}{3}-\frac{2s+2t}{6}=\frac{t}{3}
  \end{equation}
  for all points $Q\neq P$.  Since $\cM(P),\cM(Q)\in\N$ we set $x:=\cM(P)=s-\tfrac{2t}{3}$ and $y:=\cM(Q)=\tfrac{t}{3}$, so that $s=x+2y$ and $t=3y$. With this we have $n=x+6y$ and 
  $s-t=x-y$.
  \begin{itemize} 
    \item $\cM(P)=0$: $x=0$, so that we can assume $y=1$, which implies $s=2$, $t=3$, $\gamma=1$, $n=6$, $t-s=1$, and $\cM=\chi_E-\chi_P$ for some point $P$. 
    \item $\cM(Q)=0$: $y=0$, so that we can assume $x=1$, which implies $s=1$, $t=0$, $\gamma=1$, $n=1$, $s-t=1$, and $\cM=\chi_P$ for some point $P$.
  \end{itemize}
\end{Example}
So, also Example~\ref{ex_three_lines_with_common_intersection} can be explained by the subspace construction in Proposition~\ref{prop_subspace}.

\begin{Example}
  \label{ex_three_lines_without_common_intersection}
  Let $\cM$ be a multiset of points in $\PG(2,2)$ uniquely characterized by $\cM(L_1)=\cM(L_2)=\cM(L_3)=s\in\N$ for three different lines $L_1,L_2,L_3$ without a common 
  intersection point,i.e.\ $L_1\cap L_2\cap L_3=\emptyset$, and $\cM(L')=t\in\N$ for all other lines. For each point $P$ that is contained on exactly two lines $L_i$ we have
  \begin{equation}
    \cM(P) = \frac{2s+t}{3}-\frac{s+3t}{6}=\frac{3s-t}{6},
  \end{equation} 
  for each point $Q$ that is contained on exactly one line $L_i$ we have
  \begin{equation}
    \cM(Q) = \frac{s+2t}{3}-\frac{2s+2t}{6}=\frac{t}{3},
  \end{equation}
  and for the unique point $R$ that is contained on none of the lines $L_i$ we have
  \begin{equation}
    \cM(R) = \frac{3t}{3}-\frac{3s+t}{6}=\frac{5t-3s}{6}.
  \end{equation}
  Since $\cM(P),\cM(Q)\in\N$ we set $x:=\cM(P)=\tfrac{3s-t}{6}$ and $y:=\cM(Q)=\tfrac{t}{3}$, so that $s=2x+y$ and $t=3y$. With this we have $n=2x+5y$ and 
  $s-t=2(x-y)$.
  \begin{itemize} 
    \item $\cM(P)=0$: $x=0$, so that we can assume $y=1$, which implies $s=1$, $t=3$, $t-s=2$, $\gamma=2$, $n=5$, and $\cM=\chi_L+2\chi_P$ for some line $L$ and some point $P\notin L$. 
    \item $\cM(Q)=0$: $y=0$, so that $x=0$ and $\cM$ is the empty multiset of points.
    \item $\cM(R)=0$: $x=2y$, so that we can assume $y=1$, which implies $x=2$, $s=4$, $t=6$, $t-s=2$, $\gamma=2$, $n=9$, and the $2$-complement of $\cM$ equals $\cM=\chi_L+2\chi_P$ for 
                      some line $L$ and some point $P\notin L$, see the case $\cM(P)=0$.    
  \end{itemize}
\end{Example}
So, also Example~\ref{ex_three_lines_without_common_intersection} can be explained by the construction in Proposition~\ref{prop_two_subspaces}. 

\medskip

In Examples~\ref{ex_one_line}--\ref{ex_three_lines_without_common_intersection} we have considered all cases of $1\le \#\overline{\cM}\le 3$ up to symmetry. The cases 
$\#\overline{\cM}\in\{0,7\}$ give one-character multisets. By considering the complement $\cM'=\chi_V-\overline{\cM}$ we see that examples for $4\le \#\overline{\cM}\le 6$ 
do not give something new. Since the dimension of the ambient space is odd, we cannot apply the construction in Proposition~\ref{prop_partial_parallelism}.

\medskip

Now let us consider the general case. Given the set $\overline{\cM}$ of hyperplanes with multiplicity $s$ we get an explicit expression for the multiplicity $\cM(P)$ of
every point $P\in\cP$ depending on the two unknown hyperplane multiplicities $s$ and $t$. 
\begin{Lemma}
  Let $\overline{\cM}$ be a set of points in $\PG(k-1,q)$, where $k\ge 3$, and $\cM$ be a multiset of points in $\PG(k-1,q)$ such that $\cM(H)=s$ if $\overline{\cM}(H^\perp)=1$ and $\cM(H)=t$ if 
  $\overline{\cM}(H^\perp)=0$ for all hyperplanes $H\in\cH$. Denoting the number of hyperplanes $H\ni P$ with $\cM(H)=s$ by $\varphi(P)$ and setting $r:=\#\overline{\cM}$, 
  $\Delta:=s-t\in\Z$, we have 
  \begin{equation}
    \cM(P)= \frac{t+\Delta\cdot\varphi(P)}{[k-1]_q}-\frac{\Delta}{q^{k-2}}\cdot \frac{[k-2]_q}{[k-1]_q}\cdot \left(r-\varphi(P)\right). \label{eq_t_Delta}
  \end{equation}
\end{Lemma}
\begin{Proof}
  Counting gives that $[k-1]_q-\varphi(P)$ hyperplanes through $P$ have multiplicity $t$, from the $q^{k-1}$ hyperplanes not containing $P$ exactly $r-\varphi(P)$ have multiplicity $\cM(H)=s$ 
  and $q^{k-1}-r+\varphi(P)$ have multiplicity $\cM(H)=t$. With this we can use Lemma~\ref{lemma_reconstruct_points_from_hyperplanes} to compute
  \begin{eqnarray*}
    \cM(P) &=& \sum_{H\in \cH\,:\, P\in H} \frac{1}{[k-1]_q}\cdot \cM(H)+ \sum_{H\in \cH\,:\, P\notin H} \frac{1}{q^{k-1}}\cdot\left(\frac{1}{[k-1]_q}-1\right)\cdot \cM(H) \\ 
    &=& \sum_{H\in \cH\,:\, P\in H} \frac{1}{[k-1]_q}\cdot \cM(H)- \sum_{H\in \cH\,:\, P\notin H} \frac{1}{q^{k-1}}\cdot \frac{q[k-2]_q}{[k-1]_q}\cdot \cM(H) \\  
    &=& t+\frac{\Delta}{[k-1]_q}\cdot \varphi(P)\,\,-\,\, \, \frac{q[k-2]_q}{[k-1]_q}\cdot t -\frac{\Delta}{q^{k-1}}\cdot \frac{q[k-2]_q}{[k-1]_q}\cdot \left(r-\varphi(P)\right) \\ 
    &=& \frac{t+\Delta\cdot\varphi(P)}{[k-1]_q}-\frac{\Delta}{q^{k-2}}\cdot \frac{[k-2]_q}{[k-1]_q}\cdot \left(r-\varphi(P)\right).  
  \end{eqnarray*}  
\end{Proof}
Note that $\varphi(P)=\overline{\cM}(P^\perp)$ for all $P\in\cP$.

\begin{Lemma}
  \label{lemma_aux_2}
  Let $\overline{\cM}$ be a set of points in $\PG(k-1,q)$, where $k\ge 3$, and $\cM$ be a multiset of points in $\PG(k-1,q)$ such that $\cM(H)=s$ if $\overline{\cM}(H^\perp)=1$ and $\cM(H)=t$ if 
  $\overline{\cM}(H^\perp)=0$ for all hyperplanes $H\in\cH$. Denote the number of hyperplanes $H\ni P$ with $\cM(H)=s$ by $\varphi(P)$ and uniquely choose $m\in\N$, $\cI\subseteq\N$ with $0\in\cI$ such that 
  $\left\{\varphi(P)\,:\,P\in\cP\right\}=\left\{m+i\,:\,i\in\cI\right\}$. If $s> t$ and there exists a point $Q\in\cP$ with $\cM(Q)=0$, then we have 
  \begin{equation}
    t= \frac{\Delta}{q^{k-2}}\cdot [k-2]_q \cdot \left(r-m\right) -\Delta\cdot m\label{eq_t} 
  \end{equation}
  and
  \begin{equation}
    \cM(P)=\frac{\Delta\cdot i}{q^{k-2}}\label{eq_point_mult_Delta_i}
  \end{equation}
  for all points $P\in \cP$ where $i:=\varphi(P)-m$, $r:=\#\overline{\cM}$, and $\Delta:=s-t\in\N_{\ge 1}$. If $\cM$ is non-repetitive, then $\Delta$ divides $q^{k-2}$.
\end{Lemma}  
\begin{Proof}
  Using $\Delta> 0$ we observe that the expression for $\cM(P)$ in Equation~(\ref{eq_t_Delta}) is increasing in $\varphi(P)$. So, we need to choose a point $Q\in\cP$ which minimizes 
  $\varphi(Q)$ to normalize using $\cM(Q)=0$ since otherwise we will obtain points with negative multiplicity. So, choosing a point $Q\in\cP$ with $\varphi(Q)=m$ we require
  $$
    0=\cM(Q)=\frac{t+\Delta\cdot m}{[k-1]_q}-\frac{\Delta}{q^{k-1}}\cdot \frac{q[k-2]_q}{[k-1]_q}\cdot \left(r-m\right),
  $$
  which yields Equation~(\ref{eq_t}). Using $i:=\varphi(P)-m$ and the expression for $t$ we compute
  \begin{eqnarray*}
    \cM(P)&=& \frac{t+\Delta\cdot(m+i)}{[k-1]_q}-\frac{\Delta}{q^{k-2}}\cdot \frac{[k-2]_q}{[k-1]_q}\cdot \left(r-m-i\right) \\ 
    &=&  \frac{\Delta}{q^{k-2}}\cdot \frac{[k-2]_q}{[k-1]_q} \cdot \left(r-m\right) -\frac{\Delta\cdot m}{[k-1]_q}+ \frac{\Delta\cdot(m+i)}{[k-1]_q}-\frac{\Delta}{q^{k-2}}\cdot \frac{[k-2]_q}{[k-1]_q}\cdot \left(r-m-i\right) \\ 
    &=& \frac{\Delta\cdot i}{[k-1]_q}+\frac{\Delta\cdot i}{q^{k-2}}\cdot \frac{[k-2]_q}{[k-1]_q} =\frac{\Delta\cdot i}{q^{k-2}}  
  \end{eqnarray*}
  for all $P\in\cP$. Note that if $f>1$ is a divisor of $\Delta$ that is coprime to $q$, then all point multiplicities of $\cM$ are divisible by $f$. If $\Delta=q^{k-2}\cdot f$ for an integer $f>1$, 
  then all point multiplicities of $\cM$ are divisible by $f$. Thus, we have that $\Delta$ divides $q^{k-2}$.
\end{Proof}
Note that $\cI=\left\{\overline{\cM}(H)-\overline{\cM}(H')\,:\,H\in\cH\right\}$, where $H'\in\cH$ is a minimizer of $\overline{\cM}(H)$.

\begin{Lemma}
  \label{lemma_aux_3}
  Let $\overline{\cM}$ be a set of points in $\PG(k-1,q)$, where $k\ge 3$ and $\cM$ be a multiset of points in $\PG(k-1,q)$ such that $\cM(H)=s$ if $\overline{\cM}(H^\perp)=1$ and $\cM(H)=t$ if 
  $\overline{\cM}(H^\perp)=0$ for all hyperplanes $H\in\cH$. 
  Using the notation from Lemma~\ref{lemma_aux_2} we set 
  \begin{eqnarray}
    g      &=& gcd\!\left(\{i\in\cI\}\cup\{q^{k-2}\}\right),\\ 
    \Delta_0 &=& q^{k-2}/g,\\ 
    t_0 &=& \frac{1}{g}\cdot [k-2]_q \cdot \left(r-m\right) -\Delta_0\cdot m,\text{ and}\\ 
    s_0 &=& t+\Delta_0. 
  \end{eqnarray}
  If $s>t$, then we have  
  $$
    \mathbb{L}(\overline{\cM})=\left\langle(s_0,t_0)\right\rangle_{\N}+\left\langle([k-1]_q,[k-1]_q)\right\rangle_{\N}. 
  $$  
\end{Lemma}
\begin{Proof}
  Setting $\mu=\mu(\cM)\in\N$ we have that 
  $\cM':=\cM-\mu\cdot\chi_V$ is a two-character multiset corresponding to $(s',t'):=(s-\mu[k-1]_q,t-\mu[k-1]_q)\in \mathbb{L}(\overline{\cM})$ and there exists a point $Q\in\cP$ with $\cM'(Q)=0$. 
  Clearly, we have $(s',t')\in\N^2$ and $s'>t'$. From Lemma~\ref{lemma_aux_2} we conclude the existence of an integer $\Delta'\in\N_{\ge 1}$ such that 
  $t'= \frac{\Delta'}{q^{k-2}}\cdot [k-2]_q \cdot \left(r-m\right) -\Delta'\cdot m$, $s'=t'+\Delta'$, and $\cM'(P)=\frac{\Delta'\cdot i}{q^{k-2}}$ for all $P\in\cP$. Since $\cM'(P)\in \N$ for 
  all $P\in\cP$ we have that $q^{k-2}$ divides $\Delta'\cdot g$, so that $\Delta_0\in\N$ divides $\Delta'$. For $f:=\Delta'/\Delta_0\in\N_{\ge 1}$ we observe that 
  $\cM'(P)$ is divisible by $f$ and we set $\cM'':=\tfrac{1}{f}\cdot\cM'$. With this, we can check that $\cM''$ is a two-character multiset corresponding to $(s_0,t_0)\in \mathbb{L}(\overline{\cM})$.    
\end{Proof}  
Note that it is not necessary to explicitly check $t_0\in\N$ since $\cM''(P)\in\N$ is sufficient to this end.

Before we consider the problem whether $\mathbb{L}(\overline{\cM})\subseteq\N^2$ contains an element $(s,t)$ with $s>t$ we treat the so far excluded case $k=2$ separately.
\begin{Lemma}
  \label{lemma_characterization_dim_2}
  Let $\overline{\cM}$ be a set of points in $\PG(1,q)$. Then, we have
  $$
    \mathbb{L}(\overline{\cM})=\left\langle(s_0,0)\right\rangle_{\N}+\left\langle(q+1,q+1)\right\rangle_{\N}, 
  $$ 
  where $s_0=0$ if $\#\overline{\cM}\in\{0,q+1\}$ and $s_0=1$ otherwise.
\end{Lemma}
\begin{Proof}
  If $\#\overline{\cM}\in\{0,q+1\}$, then a two-character multiset $\cM$ corresponding to $(s,t)\in \overline{\cM}$ actually is a one-character multiset and there exist some integer $x\in\N$ such that 
  $\cM=x\cdot\chi_v$.
  
  Otherwise we observe that in $\PG(1,q)$ points and hyperplanes coincide and the image of $\overline{\cM}$ is $\{0,1\}$. Note that we have $\cM=t\cdot \chi_V+\sum_{P\in\cP}(s-t)\cdot\overline{\cM}(P)\cdot \chi_P$ for 
  each two-character multiset $\cM$ corresponding to $(s,t)\in \mathbb{L}(\overline{\cM})$ by definition. We can easily check $(s,t)\in \left\langle(1,0)\right\rangle_{\N}+\left\langle(q+1,q+1)\right\rangle_{\N}$.  
    The proof is completed by choosing $s=1$ and $t=0$ in our representation of $\cM$.
\end{Proof}

\begin{Theorem}
  \label{main_theorem}
  Let $\overline{\cM}$ be a set of points in $\PG(k-1,q)$ with $\#\overline{\cM}\notin\left\{0,[k]_q\right\}$, where $k\ge 2$. Then
  \begin{equation}
    \cM:=\sum_{H\in \cH} \overline{\cM}(H^\perp)\cdot \chi_H
  \end{equation}
  is a two-character multiset corresponding to $(s,t)\in \mathbb{L}(\overline{\cM})$ with $n=\#\cM=r[k-1]_q$, where $r:=\#\overline{\cM}$, $t=r[k-2]_q$, and $s=r[k-2]_q+q^{k-2}$. 
  Setting $\mu:=\mu(\cM)$ and $g:=\operatorname{gcd}(\left\{\cM(P)-\mu\,:\,P\in\cP\right\})$ the function
  \begin{equation}
    \label{eq_smallest_two_character}
    \cM':=\frac{1}{g}\cdot \left( -\mu\cdot \chi_V+ \sum_{H\in \cH} \overline{\cM}(H^\perp)\cdot \chi_H\right)= \frac{1}{g}\cdot \left( \cM-\mu\cdot \chi_V\right)
  \end{equation} 
  is a two-character multiset corresponding to $(s_0,t_0)\in \mathbb{L}(\overline{\cM})$ with $n'=\#\cM'=\tfrac{1}{q}\cdot\left(r[k-1]_q-\mu[k]_q\right)$, where $r:=\#\overline{\cM}$, 
  $t_0=\tfrac{1}{g}\cdot\left(r[k-2]_q-\mu[k-1]_q\right)$, and $s_0=\tfrac{1}{g}\cdot\left(r[k-2]_q-\mu[k-1]_q+q^{k-2}\right)$, and $g$ divides $q^{k-2}$. Moreover, we have
  \begin{equation}
    \mathbb{L}(\overline{\cM})=\left\langle(s_0,t_0)\right\rangle_{\N}+\left\langle([k-1]_q,[k-1]_q)\right\rangle_{\N},
  \end{equation}  
\end{Theorem}
\begin{Proof}
  We can easily check $\cM(H)=r[k-2]_q=t$ if $\cM(H^\perp)=0$ and $\cM(H)=r[k-2]_q+q^{k-2}=s$ if $\cM(H^\perp)=1$ for all $H\in\cH$ as well as $\#\cM=r[k-1]_q$ directly from the 
  definition of $\cM$. Using Lemma~\ref{lemma_minus_ambient_space} and Lemma~\ref{lemma_gcd} we conclude that $\cM'$ is a two-character multiset with the stated parameters. 
 
  For $k=2$ Lemma~\ref{lemma_characterization_dim_2} our last statement. For $k\ge 3$ we can apply Lemma~\ref{lemma_aux_2} to conclude $g=\operatorname{gcd}(\{i\in\cI\})$
  and use the proof of Lemma~\ref{lemma_aux_3} to conclude our last statement. Since $s,t\in\N$ and $s>t$ we have that $g$ divides $g(s-t)=q^{k-2}$.  
\end{Proof}

Using the notation from Lemma~\ref{lemma_aux_2} applied to to the multiset of points $\cM-\mu\cdot\chi_V$ from Theorem~\ref{main_theorem} we observe $\#\cI\ge 2$ for   
$\#\overline{\cM}\notin\left\{0,[k]_q\right\}$. Using the fact that $g:=\operatorname{gcd}(\left\{\cM(P)-\mu\,:\,P\in\cP\right\})$, that $g$ divides $q^{k-2}$, and 
Equation~(\ref{eq_point_mult_Delta_i}) we conclude
\begin{equation}
  g= gcd\!\left(\{i\in\cI\}\right)=gcd\!\left(\left\{\overline{\cM}(H)-\overline{\cM}(H')\,:\,H\in\cH\right\}\right),\label{eq_q_hyperplane_mult_difference}   
\end{equation}    
where $H'\in\cH$ is a minimizer of $\overline{\cM}(H)$.

Using the classification of one-character multisets we conclude from Theorem~\ref{main_theorem}:
\begin{Corollary}
  Let $\overline{\cM}$ be a set of points in $\PG(k-1,q)$, where $k\ge 2$. Then, there exist $\left(s_0,t_0\right)\in\N^2$ such that 
  $\mathbb{L}(\overline{\cM})=\left\langle(s_0,t_0)\right\rangle_{\N}+\left\langle([k-1]_q,[k-1]_q)\right\rangle_{\N}$.
\end{Corollary}

\begin{Theorem}
  \label{main_thm_2}
  Let $\widetilde{\cM}$ be a two-character multiset in $\PG(k-1,q)$, where $k\ge 2$. Then, there exist unique $u,v\in\N$ such that $\widetilde{\cM}=u\cdot \cM'+v\cdot\chi_V$, 
  where $\cM'$ is given by Equation~(\ref{eq_smallest_two_character}). 
\end{Theorem}
\begin{Proof}
  Let $s>t$ be the two hyperplane multiplicities of $\widetilde{\cM}$. With this define $\overline{\cM}$ such that $\overline{\cM}(H^\perp)=1$ if $\widetilde{\cM}(H)=s$ and 
  $\overline{\cM}(H^\perp)=0$ if $\widetilde{\cM}(H)=t$ for all $H\in\cH$. So, $(s,t)\in \mathbb{L}(\overline{\cM})$ and Theorem~\ref{main_theorem} yields the 
  existence of $u,v\in\N$ with $(s,t)=u\cdot(s_0,t_0)+v\cdot([k-1]_q,[k-1]_q)$, where $s_0$, $t_0$ are as in Theorem~\ref{main_theorem}. From Lemma~\ref{lemma_reconstruct_points_from_hyperplanes} 
  we then conclude $\widetilde{\cM}=u\cdot \cM'+v\cdot\chi_V$. Note that $\mu(\cM')$ and $\mu(\chi_V)=1$ imply $\mu(\widetilde{M})=v$, so that $u$ can be computed from 
  $\gamma(\widetilde{\cM})=u\cdot\gamma(\cM')+v$.   
\end{Proof}
Note that for a one-character multiset $\widetilde{\cM}$ there exists a unique $v\in\N$ such that $\widetilde{\cM}=v\cdot\chi_V$. Given a set of points $\overline{\cM}$ we call $\cM'$ the 
canonical representant of the set of two-character multisets $\cM$ corresponding to $(s,t)\in\mathbb{L}(\overline{\cM})$. If $\cM=\cM'$ we just say that $\cM$ is the canonical two-character multiset. 

\begin{Theorem}
  \label{main_thm_3}
  Let $w_1<w_2$ be the non-zero weights of a non-repetitive $[n,k]_q$ two-weight code $C$ without full support. Then, there exist integers $f$ and $u$ such that 
  $w_1=up^f$ and $w_2=(u+1)p^f$, where $p$ is the characteristic of the underlying field $\F_q$.
\end{Theorem}
\begin{Proof}
  Let $\cM$ be the two-character multiset in $\PG(k-1,q)$ corresponding to $C$. Choose unique $u,v\in\N$ such that $\cM=u\cdot\cM'+v\cdot\chi_V$ as in Theorem~\ref{main_thm_2}. Since $C$ does not 
  have full support, we   have $v=0$ and since $C$ is non-repetitive we have $u=1$. With this we can use Theorem~\ref{main_theorem} to compute 
  \begin{equation}
    w_1=n-s=\tfrac{1}{g}\cdot\left(r\cdot q^{k-2}-\mu\cdot q^{k-1}-q^{k-2}\right)=(r-q\mu-1)\cdot p^f
  \end{equation}
  and
  \begin{equation}
    w_2=n-t=\tfrac{1}{g}\cdot\left(r\cdot q^{k-2}-\mu\cdot q^{k-1}\right)=(r-q\mu)\cdot p^f,
  \end{equation}      
  where $f$ is chosen such that $\tfrac{q^{k-2}}{g}=p^f$. I.e., we can choose $u=r-q\mu-1$.
\end{Proof}

We have seen in Equation~(\ref{eq_q_hyperplane_mult_difference}) that we can compute the parameter $g$ directly from the set of points $\overline{\cM}$. If we additionally assume that 
$\overline{\cM}$ is spanning, then we can consider the corresponding projective $[n,k]_q$-code $\overline{C}$, where $n=\#\overline{\cM}$.\footnote{If $\overline{\cM}$ is not spanning 
then we can consider the lowerdimensional subspace spanned by $\operatorname{supp}(\overline{\cM})$.} Note that we have $\overline{\cM}(H)\equiv m\pmod g$ for all $H\in\cH$ and that $g$ 
is maximal with this property. If $m\equiv n\pmod g$, then $g$ would simply be the maximal divisibility constant of the weights of $\overline{C}$. From \cite[Theorem 7]{honold2018partial} 
or \cite[Theorem 3]{ward1999introduction} we can conclude $m\equiv n\pmod g$. Thus, we have
\begin{equation}
  g=\operatorname{gcd}\!\left(\left\{\operatorname{wt}(c)\,:\,c\in\overline{C}\right\}\right).\label{eq_divisible}
\end{equation} 
The argument may also be based on the following lemma (using the fact that $\overline{C}$ is projective):
\begin{Lemma}
  Let $C$ be an $[n,k]_q$-code of full length such that we have $\operatorname{wt}(c)\equiv m\pmod\Delta$ for all non-zero codewords $c\in C$. If $\Delta$ is a power 
  of the characteristic of the underlying field $\F_q$, then we have $m\equiv 0\pmod {\min\{\Delta,q\}}$. Moreover, if additionally $q$ divides $\Delta$ and $k\ge 2$, then 
  the non-zero weights in each residual code are congruent to $m/q$ modulo $\Delta/q$.
\end{Lemma}
\begin{Proof}
  Let $\cM$ be the multiset of points in $\PG(k-1,q)$ corresponding to $C$. For each hyperplane $H$ we have $n-\cM(H)\equiv m\pmod\Delta$, which is equivalent to 
  $\cM(H)\equiv n-m \pmod\Delta$. The weight of a non-zero codeword in a residual code is given by a subspace $K$ of codimension $2$ and a hyperplane $H$ with $K\le H$. With this, 
  the weight is given by $\cM(H)-\cM(K)\equiv n-m-\cM(K)\pmod\Delta$. Counting the hyperplane multiplicities of the $q+1$ hyperplanes that contain $K$ yields
  \begin{equation}
    \sum_{H\in\cH\,:\, K\le H} \cM(H)=\#\cM+q\cdot\cM(K)=\#\cM+q\cdot\cM(K)
  \end{equation}
  and
  \begin{equation}
    \sum_{H\in\cH\,:\, K\le H} \cM(H) \equiv(q+1)(n-m)\pmod \Delta,
  \end{equation}  
  so that
  \begin{equation}
    m\equiv q\cdot\left(n-m-\cM(K)\right)\pmod \Delta.
  \end{equation}
\end{Proof}
Given Equation~(\ref{eq_divisible}) we might be interested in projective divisible codes (with a large divisibility constant). For enumerations for the binary case we refer to 
\cite{projective_divisible_binary_codes} and for a more general survey we refer to e.g.\ \cite{kurz2021divisible}. Note that the only point sets $\cM$ in $\PG(k-1,q)$ that 
are $q^{k-1}$-divisible are given by $\#\cM\in\left\{0,[k]_q\right\}$, i.e., the empty and the full set. All other point sets are at most $q^{k-2}$-divisible, as implied by     
Theorem~\ref{main_theorem}.

\section{Enumeration of two-character multisets in $\PG(k-1,q)$ for small parameters}
\label{sec_enumeration}
Since all two-character multisets in $\PG(1,q)$ can be parameterized as $\cM=b\cdot \chi_V+\sum_{P\in\cP}(a-b)\cdot\overline{\cM}(P)\cdot \chi_P$ for integers $a>b\ge 0$ and 
a set of points $\overline{\cM}$ in $\PG(k-1,q)$, see Lemma~\ref{lemma_characterization_dim_2} and its proof, we assume $k\ge 3$ in the following. Due to Theorem~\ref{main_thm_2} every 
two-character multiset in $\PG(k-1,q)$ can be written as $u\cdot\cM'+v\cdot\chi_V$, where $u,v\in\N$ and $\cM'$ is characterized in Theorem~\ref{main_theorem}. So, we further restrict 
out considerations on canonical two-character multisets where we have $u=1$ and $v=0$. For $k=2$ all canonical two-character multisets in $\PG(k-1,q)$ are indeed sets of points and given by 
the construction in Proposition~\ref{prop_partial_parallelism} (with $r=1$).

It can be easily checked that two isomorphic sets of points in $\PG(k-1,q)$ yield isomorphic canonical two-character multisets $\cM'$. So, for the full enumeration 
of canonical two-character multisets in $\PG(k-1,q)$ we just need to loop over all non-isomorphic sets of points $\overline{\cM}$ in $\PG(k-1,q)$ and use Theorem~\ref{main_theorem} 
to determine $\cM$, $\cM'$, and their parameters. We remark that the numbers of non-isomorphic projective codes per length, dimension, and field size are e.g.\ listed in 
\cite[Tables 6.10--6.12]{betten2006error} (for small parameters). For the binary case and dimensions at most six some additional data can be found in \cite{bouyukliev2006binary}. Here 
we utilize the software package \textsc{LinCode} \cite{bouyukliev2021computer} to enumerate these codes.  

\begin{table}[htp]
  \begin{center}
    \begin{tabular}{|rrrrrrr|rrrr|}
      \hline
      $g$ & $\mu$ & $r$ & $n$ & $\gamma$ & $s$ & $t$ & $s_0$ & $t_0$ & $n'$ & $\gamma'$ \\
      \hline
2 & 1 & 3 & 9 & 3 & 5 & 3 & 1 & 0 & 1 & 1 \\
1 & 0 & 1 & 3 & 1 & 3 & 1 & 3 & 1 & 3 & 1 \\
1 & 2 & 6 & 18 & 3 & 8 & 6 & 2 & 0 & 4 & 1 \\
2 & 0 & 4 & 12 & 2 & 6 & 4 & 3 & 2 & 6 & 1 \\
\hline
1 & 1 & 4 & 12 & 3 & 6 & 4 & 3 & 1 & 5 & 2 \\
1 & 0 & 2 & 6 & 2 & 4 & 2 & 4 & 2 & 6 & 2 \\
1 & 1 & 5 & 15 & 3 & 7 & 5 & 4 & 2 & 8 & 2 \\
1 & 0 & 3 & 9 & 2 & 5 & 3 & 5 & 3 & 9 & 2 \\
\hline      
    \end{tabular}  
    \caption{Feasible parameters for canonical two-character multisets in $\PG(2,2)$.}
    \label{table_feasible_parameters_pg_2_2}
  \end{center}
\end{table}

In Table~\ref{table_feasible_parameters_pg_2_2} and in Table~\ref{table_feasible_parameters_pg_3_2} we list the feasible parameters for canonical two-character multisets in $\PG(2,2)$ and 
in $\PG(3,2)$, respectively, where $n':=\#\cM'$ and $\gamma':=\gamma(\cM')$. For $\PG(2,2)$ we can also state more direct constructions:
\begin{itemize}
  \item $(n',s_0,t_0,\gamma')=(1,1,0,1)$: characteristic function of a point (not spanning)
  \item $(n',s_0,t_0,\gamma')=(3,3,1,1)$: characteristic function of a line (not spanning)
  \item $(n',s_0,t_0,\gamma')=(4,2,0,1)$: complement of the characteristic function of a line
  \item $(n',s_0,t_0,\gamma')=(6,3,2,1)$: complement of the characteristic function of a point
  \item $(n',s_0,t_0,\gamma')=(5,3,1,2)$: $\chi_L+2\chi_P$ for a line $L$ and a point $P$ with $P\notin L$
  \item $(n',s_0,t_0,\gamma')=(6,4,2,2)$: $\chi_L+\chi_{L}'$ for two different lines $L$ and $L'$
  \item $(n',s_0,t_0,\gamma')=(8,4,2,2)$:  $\chi_V-\chi_L-\chi_{L}'$ for two different lines $L$ and $L'$ 
  \item $(n',s_0,t_0,\gamma')=(9,5,3,2)$: $2\chi_V-\chi_L+-\chi_P$ for a line $L$ and a point $P$ with $P\notin L$
\end{itemize}  

\begin{table}[htp]
  \begin{center}
    \begin{tabular}{|rrrrrrr|rrrr|}
      \hline
      $g$ & $\mu$ & $r$ & $n$ & $\gamma$ & $s$ & $t$ & $s_0$ & $t_0$ & $n'$ & $\gamma'$ \\
      \hline
4 & 3 & 7 & 49 & 7 & 25 & 21 & 1 & 0 & 1 & 1 \\
2 & 1 & 3 & 21 & 3 & 13 & 9 & 3 & 1 & 3 & 1 \\
2 & 4 & 10 & 70 & 6 & 34 & 30 & 3 & 1 & 5 & 1 \\
2 & 2 & 6 & 42 & 4 & 22 & 18 & 4 & 2 & 6 & 1 \\
1 & 0 & 1 & 7 & 1 & 7 & 3 & 7 & 3 & 7 & 1 \\
1 & 6 & 14 & 98 & 7 & 46 & 42 & 4 & 0 & 8 & 1 \\
2 & 3 & 9 & 63 & 5 & 31 & 27 & 5 & 3 & 9 & 1 \\
2 & 1 & 5 & 35 & 3 & 19 & 15 & 6 & 4 & 10 & 1 \\
2 & 4 & 12 & 84 & 6 & 40 & 36 & 6 & 4 & 12 & 1 \\
4 & 0 & 8 & 56 & 4 & 28 & 24 & 7 & 6 & 14 & 1 \\
\hline
2 & 2 & 8 & 56 & 6 & 28 & 24 & 7 & 5 & 13 & 2 \\
1 & 0 & 2 & 14 & 2 & 10 & 6 & 10 & 6 & 14 & 2 \\
2 & 0 & 4 & 28 & 4 & 16 & 12 & 8 & 6 & 14 & 2 \\
1 & 5 & 13 & 91 & 7 & 43 & 39 & 8 & 4 & 16 & 2 \\
2 & 3 & 11 & 77 & 7 & 37 & 33 & 8 & 6 & 16 & 2 \\
2 & 1 & 7 & 49 & 5 & 25 & 21 & 9 & 7 & 17 & 2 \\
\hline
1 & 1 & 4 & 28 & 4 & 16 & 12 & 9 & 5 & 13 & 3 \\
1 & 4 & 11 & 77 & 7 & 37 & 33 & 9 & 5 & 17 & 3 \\
1 & 3 & 9 & 63 & 6 & 31 & 27 & 10 & 6 & 18 & 3 \\
1 & 2 & 7 & 49 & 5 & 25 & 21 & 11 & 7 & 19 & 3 \\
1 & 1 & 5 & 35 & 4 & 19 & 15 & 12 & 8 & 20 & 3 \\
1 & 0 & 3 & 21 & 3 & 13 & 9 & 13 & 9 & 21 & 3 \\
1 & 4 & 12 & 84 & 7 & 40 & 36 & 12 & 8 & 24 & 3 \\
1 & 3 & 10 & 70 & 6 & 34 & 30 & 13 & 9 & 25 & 3 \\
1 & 2 & 8 & 56 & 5 & 28 & 24 & 14 & 10 & 26 & 3 \\
1 & 1 & 6 & 42 & 4 & 22 & 18 & 15 & 11 & 27 & 3 \\
1 & 0 & 4 & 28 & 3 & 16 & 12 & 16 & 12 & 28 & 3 \\
1 & 3 & 11 & 77 & 6 & 37 & 33 & 16 & 12 & 32 & 3 \\
\hline
1 & 3 & 8 & 56 & 7 & 28 & 24 & 7 & 3 & 11 & 4 \\
1 & 2 & 6 & 42 & 6 & 22 & 18 & 8 & 4 & 12 & 4 \\
1 & 3 & 9 & 63 & 7 & 31 & 27 & 10 & 6 & 18 & 4 \\
1 & 2 & 7 & 49 & 6 & 25 & 21 & 11 & 7 & 19 & 4 \\
1 & 1 & 5 & 35 & 5 & 19 & 15 & 12 & 8 & 20 & 4 \\
1 & 3 & 10 & 70 & 7 & 34 & 30 & 13 & 9 & 25 & 4 \\
1 & 2 & 8 & 56 & 6 & 28 & 24 & 14 & 10 & 26 & 4 \\
1 & 1 & 6 & 42 & 5 & 22 & 18 & 15 & 11 & 27 & 4 \\
1 & 2 & 9 & 63 & 6 & 31 & 27 & 17 & 13 & 33 & 4 \\
1 & 1 & 7 & 49 & 5 & 25 & 21 & 18 & 14 & 34 & 4 \\
1 & 0 & 5 & 35 & 4 & 19 & 15 & 19 & 15 & 35 & 4 \\
1 & 2 & 10 & 70 & 6 & 34 & 30 & 20 & 16 & 40 & 4 \\
1 & 1 & 8 & 56 & 5 & 28 & 24 & 21 & 17 & 41 & 4 \\
1 & 0 & 6 & 42 & 4 & 22 & 18 & 22 & 18 & 42 & 4 \\
1 & 1 & 9 & 63 & 5 & 31 & 27 & 24 & 20 & 48 & 4 \\
1 & 0 & 7 & 49 & 4 & 25 & 21 & 25 & 21 & 49 & 4 \\
      \hline
    \end{tabular}  
    \caption{Feasible parameters for canonical two-character multisets in $\PG(3,2)$.}
    \label{table_feasible_parameters_pg_3_2}
  \end{center}
\end{table}

Of course, also for $\PG(3,2)$ some of the examples have nicer descriptions:
\begin{itemize}
  \item $(n',s_0,t_0,\gamma')=(1,1,0,1)$: characteristic function of a point (not spanning)
  \item $(n',s_0,t_0,\gamma')=(3,3,1,1)$: characteristic function of a line (not spanning)
  \item $(n',s_0,t_0,\gamma')=(7,7,3,1)$: characteristic function of a plane (not spanning)
  \item $(n',s_0,t_0,\gamma')=(5,3,1,1)$: projective base; spanning projective $2$-weight code
  \item $(n',s_0,t_0,\gamma')=(6,4,2,1)$: characteristic function of two disjoint lines; spanning projective $2$-weight code
  \item $(n',s_0,t_0,\gamma')=(14,10,6,2)$: characteristic function of two different planes
  \item $(n',s_0,t_0,\gamma')=(21,13,9,3)$: characteristic function of three planes intersecting in a common point but not a common line 
\end{itemize}  

Note that we may restrict our considerations to $r<[k]_q/2$ since if $\cM'$ is the a canonical two-character multiset for a set of points $\overline{\cM}$ with $\#\overline{\cM}=r$, 
then the complement of $\cM'$ is the the a canonical two-character multiset for a set of points which is the complement of $\overline{\cM}$ and has cardinality $[k]_q-r$. 

From the data in Table~\ref{table_feasible_parameters_pg_2_2} and Table~\ref{table_feasible_parameters_pg_3_2} we can guess the the maximum possible point 
multiplicity $\gamma(\cM')$ of $\cM'$:
\begin{Proposition}
  Let $\cM$ be a canonical two-character multiset in $\PG(k-1,q)$, where $k\ge 2$. Then, we have $\gamma(\cM)\le q^{k-2}$.
\end{Proposition}
\begin{Proof}
  Choose a suitable set $\cH'\subseteq \cH$ and $g,\nu\in\N$ such that
  $$
    \cM=\frac{1}{g}\cdot\left(\sum_{H\in\cH'} \chi_H\,-\,\mu\cdot \chi_V\right).
  $$
  Let $P\in\cP$ be a point with $\cM(P)=\gamma$ and $Q\in\cP$ be a point with $\cM(Q)=0$. With this we have $\lambda\ge \left|\left\{H\in\cH'\,:\,Q\le H\right\} \right|$. 
  Since $P$ is contained $[k-1]_q$ hyperplanes in $\cH$ and $\langle P,Q\rangle$ is contained in $[k-2]_q$ hyperplanes in $\cH$ we have 
  $\cM(P)\le q^{k-2}$.
\end{Proof}

We can easily construct an example showing that the stated upper bound is tight. To this end let $P$, $Q$ be two different points in $\PG(k-1,q)$, where $k\ge 3$,and $H'$ be an arbitrary hyperplane 
neither containing $P$ nor $Q$. With this, we choose $\cH'$ as the set of all $q^{k-2}$ hyperplanes that contain $P$ but do not contain $Q$ and additionally the hyperplane $H'$. 
For the corresponding multiset of points $\cM$ we then have $\cM(P)=q^{k-2}$ and $\cM(Q)=0$, so that $\mu(\cM)=0$. For an arbitrary point $R\in H'$ we have 
$\cM(R)=q^{k-2}-q^{k-3}+1=(q-1)q^{k-3}+1$, so that $\operatorname{gcd}(\cM(R),\cM(P))=1$ if $k\ge 4$ or $k=3$ and $q\neq 2$. For $(k,q)=(3,2)$ we have already seen examples of canonical two-character 
multisets with maximum point multiplicity $2$. 

In Table~\ref{table_feasible_parameters_pg_4_2_1} and Table~\ref{table_feasible_parameters_pg_4_2_2} we list the feasible parameters for canonical two-character multisets in $\PG(4,2)$ with 
point multiplicity at most $4$.

\begin{table}[htp]
  \begin{center}
    \begin{tabular}{|rrrrrrr|rrrr|}
      \hline
      $g$ & $\mu$ & $r$ & $n$ & $\gamma$ & $s$ & $t$ & $s_0$ & $t_0$ & $n'$ & $\gamma'$ \\
      \hline
8 & 7 & 15 & 225 & 15 & 113 & 105 & 1 & 0 & 1 & 1 \\
4 & 3 & 7 & 105 & 7 & 57 & 49 & 3 & 1 & 3 & 1 \\
2 & 1 & 3 & 45 & 3 & 29 & 21 & 7 & 3 & 7 & 1 \\
1 & 0 & 1 & 15 & 1 & 15 & 7 & 15 & 7 & 15 & 1 \\
1 & 14 & 30 & 450 & 15 & 218 & 210 & 8 & 0 & 16 & 1 \\
2 & 12 & 28 & 420 & 14 & 204 & 196 & 12 & 8 & 24 & 1 \\
4 & 8 & 24 & 360 & 12 & 176 & 168 & 14 & 12 & 28 & 1 \\
8 & 0 & 16 & 240 & 8 & 120 & 112 & 15 & 14 & 30 & 1 \\
      \hline
2 & 2 & 6 & 90 & 6 & 50 & 42 & 10 & 6 & 14 & 2 \\
2 & 5 & 13 & 195 & 9 & 99 & 91 & 12 & 8 & 20 & 2 \\
2 & 3 & 9 & 135 & 7 & 71 & 63 & 13 & 9 & 21 & 2 \\
2 & 1 & 5 & 75 & 5 & 43 & 35 & 14 & 10 & 22 & 2 \\
2 & 8 & 20 & 300 & 12 & 148 & 140 & 14 & 10 & 26 & 2 \\
2 & 6 & 16 & 240 & 10 & 120 & 112 & 15 & 11 & 27 & 2 \\
2 & 4 & 12 & 180 & 8 & 92 & 84 & 16 & 12 & 28 & 2 \\
2 & 2 & 8 & 120 & 6 & 64 & 56 & 17 & 13 & 29 & 2 \\
4 & 0 & 8 & 120 & 8 & 64 & 56 & 16 & 14 & 30 & 2 \\
2 & 0 & 4 & 60 & 4 & 36 & 28 & 18 & 14 & 30 & 2 \\
1 & 0 & 2 & 30 & 2 & 22 & 14 & 22 & 14 & 30 & 2 \\
4 & 7 & 23 & 345 & 15 & 169 & 161 & 16 & 14 & 32 & 2 \\
2 & 11 & 27 & 405 & 15 & 197 & 189 & 16 & 12 & 32 & 2 \\
1 & 13 & 29 & 435 & 15 & 211 & 203 & 16 & 8 & 32 & 2 \\
2 & 9 & 23 & 345 & 13 & 169 & 161 & 17 & 13 & 33 & 2 \\
2 & 7 & 19 & 285 & 11 & 141 & 133 & 18 & 14 & 34 & 2 \\
2 & 5 & 15 & 225 & 9 & 113 & 105 & 19 & 15 & 35 & 2 \\
2 & 3 & 11 & 165 & 7 & 85 & 77 & 20 & 16 & 36 & 2 \\
2 & 10 & 26 & 390 & 14 & 190 & 182 & 20 & 16 & 40 & 2 \\
2 & 8 & 22 & 330 & 12 & 162 & 154 & 21 & 17 & 41 & 2 \\
2 & 6 & 18 & 270 & 10 & 134 & 126 & 22 & 18 & 42 & 2 \\
2 & 9 & 25 & 375 & 13 & 183 & 175 & 24 & 20 & 48 & 2 \\
      \hline
2 & 4 & 10 & 150 & 10 & 78 & 70 & 9 & 5 & 13 & 3 \\
2 & 9 & 21 & 315 & 15 & 155 & 147 & 10 & 6 & 18 & 3 \\
2 & 3 & 9 & 135 & 9 & 71 & 63 & 13 & 9 & 21 & 3 \\
2 & 6 & 16 & 240 & 12 & 120 & 112 & 15 & 11 & 27 & 3 \\
2 & 4 & 12 & 180 & 10 & 92 & 84 & 16 & 12 & 28 & 3 \\
2 & 2 & 8 & 120 & 8 & 64 & 56 & 17 & 13 & 29 & 3 \\
1 & 1 & 4 & 60 & 4 & 36 & 28 & 21 & 13 & 29 & 3 \\
2 & 7 & 19 & 285 & 13 & 141 & 133 & 18 & 14 & 34 & 3 \\
2 & 5 & 15 & 225 & 11 & 113 & 105 & 19 & 15 & 35 & 3 \\
2 & 3 & 11 & 165 & 9 & 85 & 77 & 20 & 16 & 36 & 3 \\
2 & 1 & 7 & 105 & 7 & 57 & 49 & 21 & 17 & 37 & 3 \\
2 & 6 & 18 & 270 & 12 & 134 & 126 & 22 & 18 & 42 & 3 \\
2 & 4 & 14 & 210 & 10 & 106 & 98 & 23 & 19 & 43 & 3 \\
2 & 2 & 10 & 150 & 8 & 78 & 70 & 24 & 20 & 44 & 3 \\
1 & 0 & 3 & 45 & 3 & 29 & 21 & 29 & 21 & 45 & 3 \\
1 & 12 & 28 & 420 & 15 & 204 & 196 & 24 & 16 & 48 & 3 \\
2 & 7 & 21 & 315 & 13 & 155 & 147 & 25 & 21 & 49 & 3 \\
\hline
    \end{tabular}  
    \caption{Feasible parameters for canonical two-character multisets in $\PG(4,2)$ with $\gamma'\le 4$ -- part 1.}
    \label{table_feasible_parameters_pg_4_2_1}
  \end{center}
\end{table}

\begin{table}[htp]
  \begin{center}
    \begin{tabular}{|rrrrrrr|rrrr|}
      \hline
      $g$ & $\mu$ & $r$ & $n$ & $\gamma$ & $s$ & $t$ & $s_0$ & $t_0$ & $n'$ & $\gamma'$ \\
      \hline
2 & 5 & 17 & 255 & 11 & 127 & 119 & 26 & 22 & 50 & 3 \\
2 & 3 & 13 & 195 & 9 & 99 & 91 & 27 & 23 & 51 & 3 \\
2 & 8 & 24 & 360 & 14 & 176 & 168 & 28 & 24 & 56 & 3 \\
2 & 6 & 20 & 300 & 12 & 148 & 140 & 29 & 25 & 57 & 3 \\
2 & 4 & 16 & 240 & 10 & 120 & 112 & 30 & 26 & 58 & 3 \\
2 & 2 & 12 & 180 & 8 & 92 & 84 & 31 & 27 & 59 & 3 \\
1 & 11 & 27 & 405 & 14 & 197 & 189 & 32 & 24 & 64 & 3 \\
2 & 7 & 23 & 345 & 13 & 169 & 161 & 32 & 28 & 64 & 3 \\
2 & 5 & 19 & 285 & 11 & 141 & 133 & 33 & 29 & 65 & 3 \\
2 & 3 & 15 & 225 & 9 & 113 & 105 & 34 & 30 & 66 & 3 \\
2 & 6 & 22 & 330 & 12 & 162 & 154 & 36 & 32 & 72 & 3 \\
2 & 0 & 10 & 150 & 6 & 78 & 70 & 39 & 35 & 75 & 3 \\
2 & 5 & 21 & 315 & 11 & 155 & 147 & 40 & 36 & 80 & 3 \\
      \hline
2 & 4 & 12 & 180 & 12 & 92 & 84 & 16 & 12 & 28 & 4 \\
1 & 2 & 6 & 90 & 6 & 50 & 42 & 20 & 12 & 28 & 4 \\
2 & 3 & 11 & 165 & 11 & 85 & 77 & 20 & 16 & 36 & 4 \\
2 & 4 & 14 & 210 & 12 & 106 & 98 & 23 & 19 & 43 & 4 \\
1 & 2 & 7 & 105 & 6 & 57 & 49 & 27 & 19 & 43 & 4 \\
1 & 1 & 5 & 75 & 5 & 43 & 35 & 28 & 20 & 44 & 4 \\
1 & 11 & 26 & 390 & 15 & 190 & 182 & 25 & 17 & 49 & 4 \\
1 & 10 & 24 & 360 & 14 & 176 & 168 & 26 & 18 & 50 & 4 \\
1 & 9 & 22 & 330 & 13 & 162 & 154 & 27 & 19 & 51 & 4 \\
1 & 8 & 20 & 300 & 12 & 148 & 140 & 28 & 20 & 52 & 4 \\
1 & 7 & 18 & 270 & 11 & 134 & 126 & 29 & 21 & 53 & 4 \\
1 & 5 & 14 & 210 & 9 & 106 & 98 & 31 & 23 & 55 & 4 \\
2 & 6 & 20 & 300 & 14 & 148 & 140 & 29 & 25 & 57 & 4 \\
1 & 3 & 10 & 150 & 7 & 78 & 70 & 33 & 25 & 57 & 4 \\
2 & 4 & 16 & 240 & 12 & 120 & 112 & 30 & 26 & 58 & 4 \\
1 & 2 & 8 & 120 & 6 & 64 & 56 & 34 & 26 & 58 & 4 \\
1 & 1 & 6 & 90 & 5 & 50 & 42 & 35 & 27 & 59 & 4 \\
1 & 0 & 4 & 60 & 4 & 36 & 28 & 36 & 28 & 60 & 4 \\
1 & 11 & 27 & 405 & 15 & 197 & 189 & 32 & 24 & 64 & 4 \\
1 & 10 & 25 & 375 & 14 & 183 & 175 & 33 & 25 & 65 & 4 \\
2 & 3 & 15 & 225 & 11 & 113 & 105 & 34 & 30 & 66 & 4 \\
1 & 9 & 23 & 345 & 13 & 169 & 161 & 34 & 26 & 66 & 4 \\
2 & 1 & 11 & 165 & 9 & 85 & 77 & 35 & 31 & 67 & 4 \\
1 & 8 & 21 & 315 & 12 & 155 & 147 & 35 & 27 & 67 & 4 \\
1 & 6 & 17 & 255 & 10 & 127 & 119 & 37 & 29 & 69 & 4 \\
1 & 4 & 13 & 195 & 8 & 99 & 91 & 39 & 31 & 71 & 4 \\
1 & 3 & 11 & 165 & 7 & 85 & 77 & 40 & 32 & 72 & 4 \\
1 & 2 & 9 & 135 & 6 & 71 & 63 & 41 & 33 & 73 & 4 \\
1 & 1 & 7 & 105 & 5 & 57 & 49 & 42 & 34 & 74 & 4 \\
1 & 0 & 5 & 75 & 4 & 43 & 35 & 43 & 35 & 75 & 4 \\
1 & 10 & 26 & 390 & 14 & 190 & 182 & 40 & 32 & 80 & 4 \\
2 & 3 & 17 & 255 & 11 & 127 & 119 & 41 & 37 & 81 & 4 \\
1 & 9 & 24 & 360 & 13 & 176 & 168 & 41 & 33 & 81 & 4 \\
2 & 4 & 20 & 300 & 12 & 148 & 140 & 44 & 40 & 88 & 4 \\
2 & 3 & 19 & 285 & 11 & 141 & 133 & 48 & 44 & 96 & 4 \\
1 & 9 & 25 & 375 & 13 & 183 & 175 & 48 & 40 & 96 & 4 \\
\hline
    \end{tabular}  
    \caption{Feasible parameters for canonical two-character multisets in $\PG(4,2)$ with $\gamma'\le 4$ -- part 2.}
    \label{table_feasible_parameters_pg_4_2_2}
  \end{center}
\end{table}

\appendix

\section{Feasible parameters for binary projective two-weight codes}
In this appendix we want to utilize our parameterization for the parameters of two-character multisets in order to determine the set of feasible parameters 
of two-character sets being equivalent to projective two-weight codes. 

Starting with the sum of the characteristic functions of $r$ different hyperplanes in $\PG(k-1,q)$, where $k\ge 3$, we have 
$n=r[k-1]_q$ for the cardinality and $s=r[k-2]_q+q^{k-2}$, $t=r[k-2]_q$ for the two different hyperplane multiplicities. Our first 
condition is
\begin{equation}
  1\le r\le [k]_q-1
\end{equation}
since since there are only $[k]_q$ different hyperplanes in $\PG(k-1,q)$, $r=0$ yields the trivial multiset of points with cardinality $0$, and $r=[k]_q$ 
yields the multiset of points $[k-1]_q\cdot\chi_V$, which is a complement if the trivial multiset of points. For the canonical multiset of points we 
also need the parameters $\mu$ and $g$ to state
\begin{equation}
  n'=\frac{1}{g}\cdot\left(r[k-1]_q-\mu[k]_q\right)
\end{equation} 
for the cardinality and 
\begin{equation}
  s_0=\frac{1}{g}\cdot\left(r[k-2]_q+q^{k-2}-\mu[k-1]_q\right), t_0=\frac{1}{g}\cdot\left(r[k-2]_q-\mu[k-1]_q\right) 
\end{equation} 
for the hyperplane multiplicities. Clearly, we have the conditions
\begin{equation}
  n',s_0\in\N_{\ge_1}\text{ and }t_0\in\N
\end{equation}
as well as
\begin{equation}
  g\text{ divides } q^{k-2}
\end{equation}
and
\begin{equation}
  n'<[k]_q
\end{equation}
since there are only $[k]_q$ different points in $PG(k-1,q)$. Let us denote the number of hyperplanes with multiplicity $s_0$ by $a_s$ and the number of hyperplanes with multiplicity $t_0$ by 
$a_t$, so that $a_s+a_t=[k]_q$ and $s_0a_s-t_0a_t=n'[k-1]_q$. Thus, we have
\begin{eqnarray}
  a_s &=& \frac{g}{q^{k-2}}\cdot\left(n'[k-1]_q-t_0[k]_q\right)\\
  a_t &=& [k]_q-a_s,
\end{eqnarray}  
so that we can require
\begin{equation}
  a_s,a_t\in\N_{\ge 1}.
\end{equation}
The, so-called, standard equations for sets of points in $\PG(k-1,q)$ are completed by
\begin{equation}
  {s_0\choose 2}a_s+{t_0\choose 2}a_t = {n'\choose 2}[k-2]_q.
\end{equation}

So, in order to enumerate all possible candidates for the parameters of two-character sets $\cM'$ in $\PG(k-1,q)$ we loop over all choices of $1\le r\le [k]_q-1$, $0\le\mu\le r[k-2]_q/[k-1]_q$, 
and all divisors $g$ of $q^{k-2}$, compute $n'$, $s_0$, $t_0$, $a_s$, $a_t$, and check the mentioned conditions. We remark that the weights of the corresponding two-weight codes $c'$ are given 
by $w_1=n'-s_0$ and $w_2=n'-t_0$, so that $C'$ is $q^{k-2}/g$-divisible and projective.

In the subsequent sections we list the obtained candidates and mention corresponding constructions from the literature. It turns out that for the binary case our necessary conditions
are also sufficient for $k\le 5$. For larger dimensions we have to exclude a few cases and also leave a few cases open for $k\in\{11,12\}$. 

Before we start to state the numerical details we would like to mention two parametric series of constructions (assuming that the parameters $k$ and $q$ are clear from the context). 
For $1\le i\le k-1$ we may consider the characteristic function $\cS_i:=\chi_S$ of an $i$-space $S$ with parameters $n'=[i]_q$, $s_0=[i]_q$, and $t_0=[i-1]_q$, see 
\cite[Construction SU]{calderbank1986geometry}, as well as the complement $\chi_V-\chi_S$ with $n'=[k]_q-[i]_q$, $s_0=[k-1]_q-[i]-1_q$, and $t_0=[k-1]_q-[i]_q$. 
Whenever the dimension $k$ of the ambient space is even we can consider the characteristic function $\cP_i$ of partial $(k/2)$-spreads consisting of $i$ pairwise disjoint $(k/2)$-spaces, 
where $0< i< [k]_q/[k/2]_q=q^{k/2}+1$, see \cite[Construction SU2]{calderbank1986geometry}. The parameters are given by $n'=i[k/2]_q$, $s_0=i[k/2-1]_q+q^{k/2-1}$, and $t_0=i[k/2-1]_q$.
Note that $\cP_1$ has the same parameters as $\cS_{k/2}$, $\cP_{[k]_q-1}$ has the same parameters as the complement of $\cS_{k/2}$, and the complement of $\cP_i$ has the same parameters as 
$\cP_{q^{k/2}+1-i}$. In general we can have lots of different non-isomorphic two-character sets with the same parameters, so that different parametric constructions may also cover the same
parameters. For a third possible series $\cK_i$ we currently need to assume $q=2$ and an even dimension $k$ of the ambient space. Moreover, no explicit construction is known. However, the parameters 
$n'=i\cdot \left(q^{k/2}+1\right)$, $s_0=i\cdot \left(q^{k/2-1}+1\right)$, and $t_0=i\cdot \left(q^{k/2-1}+1\right)-q^{k/2-1}$ show a clear pattern. In principle we may choose 
$0<i<[k/2]_q$ but it will turn out that rather small and rather large values of $i$ have to be excluded. Note that the complementary parameters of $\cK_i$ are given by $\cK_{[k/2]_q-i}$. We will 
discuss explicit constructions attaining the parameters of $\cK_i$ when considering explicit dimensions $k\equiv 0\pmod 2$. Finally, we remark that for a two-character set $\cM$ the corresponding 
set $\overline{\cM}$ of hyperplanes with the maximum hyperplane multiplicity is also a two-character set, so that interchanging the roles of $\cM$ and $\overline{\cM}$ gives another general construction.   

\subsection{Feasible parameters of two-character sets in $\mathbf{\PG(3,2)}$}
In Table~\ref{table_feasible_parameters_sets_pg_3_2} we have listed the feasible parameters of two-character sets in $\PG(3,2)$. We have tried to arrange those 
parameters in several series so that patterns become visible. While parameters of the complementary set of parameters can be easily computed we will mostly state them nevertheless. 
The first block consists of the parameters for the multisets $\cS_i$ for $1\le i\le 3$ obtained from the subspace construction and their complements. 
The second block consists of the parameters for the multisets $\cP_i$ for $2\le i\le 3$ obtained from the partial spread construction. Note that $\cP_1$ has the same parameters as $\cS_2$, 
$\cP_4$ has the same parameters as the complement of $\cS_2$, and the complement of $\cP_i$ has the same parameters as $\cP_{5-i}$. The third block has the parameters of $\cK_i$ for $1\le i\le 2$.

\begin{table}[htp]
  \begin{center}
    \begin{tabular}{|rrr|rrr|rr|l|}
      \hline
      $r$ & $\mu$ & $g$ & $n'$ & $s_0$ & $t_0$ & $a_s$ & $a_t$ & construction \\
      \hline 
       7 & 3 & 4 &  1 & 1 & 0 &  7 &  8 & \cite[Construction SU]{calderbank1986geometry} \\ 
       3 & 1 & 2 &  3 & 3 & 1 &  3 & 12 & \cite[Construction SU]{calderbank1986geometry} \\
       1 & 0 & 1 &  7 & 7 & 3 &  1 & 14 & \cite[Construction SU]{calderbank1986geometry} \\
      14 & 6 & 1 &  8 & 4 & 0 & 14 &  1 & complement \\
      12 & 4 & 2 & 12 & 6 & 4 & 12 &  3 & complement \\             
       8 & 0 & 4 & 14 & 7 & 6 &  8 &  7 & complement \\ 
      \hline
       6 & 2 & 2 & 6 & 4 & 2 & 6 & 9 & \cite[Construction SU2]{calderbank1986geometry} \\
       9 & 3 & 2 & 9 & 5 & 3 & 9 & 6 & complement \\ 
      \hline
      10 & 4 & 2 &  5 & 3 & 1 & 10 & 5  & \cite[Construction RT2]{calderbank1986geometry} \\ 
       5 & 1 & 2 & 10 & 6 & 4 &  5 & 10 & complement \\  
      \hline
    \end{tabular}  
    \caption{Feasible parameters for two-character sets in $\PG(3,2)$.}
    \label{table_feasible_parameters_sets_pg_3_2}
  \end{center}
\end{table}

\subsection{Feasible parameters of two-character sets in $\mathbf{\PG(4,2)}$}
In Table~\ref{table_feasible_parameters_sets_pg_4_2} we have listed the feasible parameters of two-character sets in $\PG(4,2)$. Here we have just one block 
with attaining examples $\cS_1,\dots,\cS_4$ and their complements.

\begin{table}[htp]
  \begin{center}
    \begin{tabular}{|rrr|rrr|rr|l|}
      \hline
      $r$ & $\mu$ & $g$ & $n'$ & $s_0$ & $t_0$ & $a_s$ & $a_t$ & construction \\
      \hline 
      15 &  7 & 8 &  1 &  1 &  0 & 15 & 16 & \cite[Construction SU]{calderbank1986geometry} \\
       7 &  3 & 4 &  3 &  3 &  1 &  7 & 24 & \cite[Construction SU]{calderbank1986geometry} \\
       3 &  1 & 2 &  7 &  7 &  3 &  3 & 28 & \cite[Construction SU]{calderbank1986geometry} \\
       1 &  0 & 1 & 15 & 15 &  7 & 1  & 30 & \cite[Construction SU]{calderbank1986geometry} \\
      30 & 14 & 1 & 16 &  8 &  0 & 30 &  1 & complement \\ 
      28 & 12 & 2 & 24 & 12 &  8 & 28 &  3 & complement \\ 
      24 &  8 & 4 & 28 & 14 & 12 & 24 &  7 & complement \\
      16 &  0 & 8 & 30 & 15 & 14 & 16 & 15 & complement \\ 
      \hline
    \end{tabular}  
    \caption{Feasible parameters for two-character sets in $\PG(4,2)$.}
    \label{table_feasible_parameters_sets_pg_4_2}
  \end{center}
\end{table}

\subsection{Feasible parameters of two-character sets in $\mathbf{\PG(5,2)}$}
In Table~\ref{table_feasible_parameters_sets_pg_5_2} we have listed the feasible parameters of two-character sets in $\PG(5,2)$. The first block consists of the parameters for the 
multisets $\cS_i$ for $1\le i\le 5$ obtained from the subspace construction and their complements. The second block consists of the parameters for the multisets $\cP_i$ for $2\le i\le 7$ 
obtained from the partial spread construction. The third block has the parameters of $\cK_i$ for $2\le i\le 5$. The parameters for $\cK_1$ as well as the complementary parameters for 
$\cK_6$ cannot be attained:
\begin{Lemma}
  \label{label_exl_k6_1}
  No $4$-divisible $[9,k]_2$-code exists.
\end{Lemma}
\begin{Proof}
  Apply e.g.\ \cite[Theorem 1]{kiermaier2020lengths}.
\end{Proof}

\begin{table}[htp]
  \begin{center}
    \begin{tabular}{|rrr|rrr|rr|l|}
      \hline
      $r$ & $\mu$ & $g$ & $n'$ & $s_0$ & $t_0$ & $a_s$ & $a_t$ & construction/non-existence \\
      \hline 
      31 & 15 & 16 &  1 &  1 &  0 & 31 & 32 & \cite[Construction SU]{calderbank1986geometry} \\
      15 &  7 &  8 &  3 &  3 &  1 & 15 & 48 & \cite[Construction SU]{calderbank1986geometry} \\
       7 &  3 &  4 &  7 &  7 &  3 &  7 & 56 & \cite[Construction SU]{calderbank1986geometry} \\
       3 &  1 &  2 & 15 & 15 &  7 &  3 & 60 & \cite[Construction SU]{calderbank1986geometry} \\
       1 &  0 &  1 & 31 & 31 & 15 &  1 & 62 & \cite[Construction SU]{calderbank1986geometry} \\   
      62 & 30 &  1 & 32 & 16 &  0 & 62 &  1 & complement \\
      60 & 28 &  2 & 48 & 24 & 16 & 60 &  3 & complement \\ 
      56 & 24 &  4 & 56 & 28 & 24 & 56 &  7 & complement \\
      48 & 16 &  8 & 60 & 30 & 28 & 48 & 15 & complement \\
      32 &  0 & 16 & 62 & 31 & 30 & 32 & 31 & complement \\
      \hline  
      14 &  6 & 4 & 14 & 10 &  6 & 14 & 49 & \cite[Construction SU2]{calderbank1986geometry} \\ 
      21 &  9 & 4 & 21 & 13 &  9 & 21 & 42 & \cite[Construction SU2]{calderbank1986geometry} \\  
      28 & 12 & 4 & 28 & 16 & 12 & 28 & 35 & \cite[Construction SU2]{calderbank1986geometry} \\
      35 & 15 & 4 & 35 & 19 & 15 & 35 & 28 & complement \\  
      42 & 18 & 4 & 42 & 22 & 18 & 42 & 21 & complement \\   
      49 & 21 & 4 & 49 & 25 & 21 & 49 & 14 & complement \\
      \hline
      54 & 26 & 4 &  9 &  5 &  1 & 54 &  9 & Lemma~\ref{label_exl_k6_1} \\
      45 & 21 & 4 & 18 & 10 &  6 & 45 & 18 & \cite{bierbrauer1997family} \\ 
      36 & 16 & 4 & 27 & 15 & 11 & 36 & 27 & \cite[Construction RT2]{calderbank1986geometry} \\
      27 & 11 & 4 & 36 & 20 & 16 & 27 & 36 & complement \\
      18 &  6 & 4 & 45 & 25 & 21 & 18 & 45 & complement \\
       9 &  1 & 4 & 54 & 30 & 26 &  9 & 54 & complement does not exist \\ 
      \hline
    \end{tabular}  
    \caption{Feasible parameters for two-character sets in $\PG(5,2)$.}
    \label{table_feasible_parameters_sets_pg_5_2}
  \end{center}
\end{table}

\subsection{Feasible parameters of two-character sets in $\mathbf{\PG(6,2)}$}
In Table~\ref{table_feasible_parameters_sets_pg_6_2} we have listed the feasible parameters of two-character sets in $\PG(6,2)$. Here we have just one block 
with attaining examples $\cS_1,\dots,\cS_6$ and their complements.

\begin{table}[htp]
  \begin{center}
    \begin{tabular}{|rrr|rrr|rr|l|}
      \hline
      $r$ & $\mu$ & $g$ & $n'$ & $s_0$ & $t_0$ & $a_s$ & $a_t$ & construction \\
      \hline 
       63 & 31 & 32 &   1 &  1 &  0 &  63 &  64 & \cite[Construction SU]{calderbank1986geometry} \\
       31 & 15 & 16 &   3 &  3 &  1 &  31 &  96 & \cite[Construction SU]{calderbank1986geometry} \\ 
       15 &  7 &  8 &   7 &  7 &  3 &  15 & 112 & \cite[Construction SU]{calderbank1986geometry} \\
        7 &  3 &  4 &  15 & 15 &  7 &   7 & 120 & \cite[Construction SU]{calderbank1986geometry} \\
        3 &  1 &  2 &  31 & 31 & 15 &   3 & 124 & \cite[Construction SU]{calderbank1986geometry} \\ 
        1 &  0 &  1 &  63 & 63 & 31 &   1 & 126 & \cite[Construction SU]{calderbank1986geometry} \\
      126 & 62 &  1 &  64 & 32 &  0 & 126 &   1 & complement \\ 
      124 & 60 &  2 &  96 & 48 & 32 & 124 &   3 & complement \\ 
      120 & 56 &  4 & 112 & 56 & 48 & 120 &   7 & complement \\ 
      112 & 48 &  8 & 120 & 60 & 56 & 112 &  15 & complement \\
       96 & 32 & 16 & 124 & 62 & 60 &  96 &  31 & complement \\  
       64 &  0 & 32 & 126 & 63 & 62 &  64 &  63 & complement \\
      \hline
    \end{tabular}  
    \caption{Feasible parameters for two-character sets in $\PG(6,2)$.}
    \label{table_feasible_parameters_sets_pg_6_2}
  \end{center}
\end{table}

\subsection{Feasible parameters of two-character sets in $\mathbf{\PG(7,2)}$}
In Table~\ref{table_feasible_parameters_sets_pg_7_2} we have listed the feasible parameters of two-character sets in $\PG(7,2)$ that are not covered by $\cP_i$ for $2\le i\le 15$, $\cS_i$ for 
$1\le i\le 7$, or their complements. There remains one block corresponding to the parameters of $\cK_i$ for $1\le i\le 14$. The parameters for $\cK_1,\cK_2$ as well as the complementary parameters for 
$\cK_{14},\cK_{13}$ cannot be attained:  
\begin{Lemma}
  \label{label_exl_k8_1}
  No $8$-divisible $[17,k]_2$-code exists. Each $8$-divisible $[34,k]_2$-code has column multiplicity at least two. 
\end{Lemma}
\begin{Proof}
  For the first statement apply e.g.\ \cite[Theorem 1]{kiermaier2020lengths} and for the second see e.g.\ \cite{restricted_column_multiplicities} or apply \cite[Theorem 12]{honold2018partial}.
\end{Proof}

\begin{table}[htp]
  \begin{center}
    \begin{tabular}{|rrr|rrr|rr|l|}
      \hline
      $r$ & $\mu$ & $g$ & $n'$ & $s_0$ & $t_0$ & $a_s$ & $a_t$ & construction/non-existence \\
      \hline 
      238 & 118 & 8 &  17 &   9 &   1 & 238 &  17 & Lemma~\ref{label_exl_k8_1} \\
      221 & 109 & 8 &  34 &  18 &  10 & 221 &  34 & Lemma~\ref{label_exl_k8_1} \\
      204 & 100 & 8 &  51 &  27 &  19 & 204 &  51 & \cite[Construction CY1]{calderbank1986geometry}, \cite{kohnert2007constructing} \\ 
      187 &  91 & 8 &  68 &  36 &  28 & 187 &  68 & \cite[Construction FE1]{calderbank1986geometry}, \cite{gulliver1997new}, \cite{kohnert2007constructing} \\ 
      170 &  82 & 8 &  85 &  45 &  37 & 170 &  85 & \cite[Construction CY1]{calderbank1986geometry}, \cite{kohnert2007constructing} \\ 
      153 &  73 & 8 & 102 &  54 &  46 & 153 & 102 & \cite[Construction CY2]{calderbank1986geometry}, \cite{kohnert2007constructing} \\
      136 &  64 & 8 & 119 &  63 &  55 & 136 & 119 & \cite[Construction RT2]{calderbank1986geometry}, \cite{kohnert2007constructing} \\
      119 &  55 & 8 & 136 &  72 &  64 & 119 & 136 & complement \\ 
      102 &  46 & 8 & 153 &  81 &  73 & 102 & 153 & complement \\
       85 &  37 & 8 & 170 &  90 &  82 &  85 & 170 & complement \\
       68 &  28 & 8 & 187 &  99 &  91 &  68 & 187 & complement \\
       51 &  19 & 8 & 204 & 108 & 100 &  51 & 204 & complement \\
       34 &  10 & 8 & 221 & 117 & 109 &  34 & 221 & complement does not exist \\ 
       17 &   1 & 8 & 238 & 126 & 118 &  17 & 238 & complement does not exist \\
      \hline
    \end{tabular}  
    \caption{Feasible parameters for two-character sets in $\PG(7,2)$ excluding those for $\cP_i$, $\cS_i$, and their complements.}
    \label{table_feasible_parameters_sets_pg_7_2}
  \end{center}
\end{table}

\subsection{Feasible parameters of two-character sets in $\mathbf{\PG(8,2)}$}
In Table~\ref{table_feasible_parameters_sets_pg_8_2} we have listed the feasible parameters of two-character sets in $\PG(8,2)$ that are not covered by $\cS_i$ for 
$1\le i\le 8$, or their complements.


\begin{table}[htp]
  \begin{center}
    \begin{tabular}{|rrr|rrr|rr|l|}
      \hline
      $r$ & $\mu$ & $g$ & $n'$ & $s_0$ & $t_0$ & $a_s$ & $a_t$ & construction \\
      \hline
      441 & 217 & 8 & 196 & 100 &  84 & 441 &  70 & \cite{bierbrauer1997family} \\ 
       73 &  33 & 8 & 219 & 123 & 107 &  73 & 438 & \cite{kohnert2007constructing} \\
      438 & 214 & 8 & 292 & 148 & 132 & 438 &  73 & complement \\ 
       70 &  30 & 8 & 315 & 171 & 155 &  70 & 441 & complement \\
      \hline
      315 & 155 & 16 &  70 &  38 &  30 & 315 & 196 & \cite{bierbrauer1997family} \\ 
      219 & 107 & 16 &  73 &  41 &  33 & 219 & 292 & \cite{fiedler1998strongly}, \cite{kohnert2007constructing} \\
      292 & 132 & 16 & 438 & 222 & 214 & 292 & 219 & complement \\      
      196 &  84 & 16 & 441 & 225 & 217 & 196 & 315 & complement \\
      \hline
    \end{tabular}  
    \caption{Feasible parameters for two-character sets in $\PG(8,2)$ excluding those for $\cS_i$ and their complements.}
    \label{table_feasible_parameters_sets_pg_8_2}
  \end{center}
\end{table}

\subsection{Feasible parameters of two-character sets in $\mathbf{\PG(9,2)}$}
In Table~\ref{table_feasible_parameters_sets_pg_9_2} we have listed the feasible parameters of two-character sets in $\PG(9,2)$ that are not covered by $\cP_i$ for $2\le i\le 31$, $\cS_i$ for 
$1\le i\le 9$, or their complements. There remains one block corresponding to the parameters of $\cK_i$ for $1\le i\le 30$. The parameters for $\cK_1,\dots,\cK_5$ as well as the complementary parameters for 
$\cK_{30},\dots,\cK_{26}$ cannot be attained:  
\begin{Lemma}
  \label{label_exl_k10_1}
  No projective $16$-divisible $[n,k]_2$-code exists for $n\in\{33,66,99,132\}$. 
\end{Lemma}
\begin{Proof}
  For $n\in \{33,66,99\}$ we apply \cite[Theorem 12]{honold2018partial}. For $n=132$ we apply \cite[Corollary 6.15]{kurz2021divisible} with $t=4$.
\end{Proof} 

\begin{Lemma}
  \label{label_exl_k10_2}
  No $[165,10,80]_2$-code exists. 
\end{Lemma}
\begin{Proof}
  See \cite[Corollary 1]{bouyukliev2000some}.
\end{Proof}
The result is concluded from the non-existence of a $[69,9,32]_2$-code applied to the residual code of a codeword of weight $96$ c.f.~\cite[Section 3.4]{guritman2000restrictions}.

\begin{table}[htp]
  \begin{center}
    \begin{tabular}{|rrr|rrr|rr|l|}
      \hline
      $r$ & $\mu$ & $g$ & $n'$ & $s_0$ & $t_0$ & $a_s$ & $a_t$ & construction/non-existence \\
      \hline 
      990 & 494 & 16 &  33 &  17 &   1 & 990 &  33 & Lemma~\ref{label_exl_k10_1} \\ 
      957 & 477 & 16 &  66 &  34 &  18 & 957 &  66 & Lemma~\ref{label_exl_k10_1} \\ 
      924 & 460 & 16 &  99 &  51 &  35 & 924 &  99 & Lemma~\ref{label_exl_k10_1} \\  
      891 & 443 & 16 & 132 &  68 &  52 & 891 & 132 & Lemma~\ref{label_exl_k10_1} \\  
      858 & 426 & 16 & 165 &  85 &  69 & 858 & 165 & Lemma~\ref{label_exl_k10_2} \\  
      825 & 409 & 16 & 198 & 102 &  86 & 825 & 198 & \cite{kohnert2007constructing} \\
      792 & 392 & 16 & 231 & 119 & 103 & 792 & 231 & \cite{dissett2000combinatorial}, \cite{kohnert2007constructing} \\
      759 & 375 & 16 & 264 & 136 & 120 & 759 & 264 & \cite[Construction FE1]{calderbank1986geometry}, \cite{kohnert2007constructing} \\ 
      726 & 358 & 16 & 297 & 153 & 137 & 726 & 297 & \cite{dissett2000combinatorial}, \cite{kohnert2007constructing} \\  
      693 & 341 & 16 & 330 & 170 & 154 & 693 & 330 & \cite{dissett2000combinatorial}, \cite{kohnert2007constructing} \\
      660 & 324 & 16 & 363 & 187 & 171 & 660 & 363 & \cite{dissett2000combinatorial}, \cite{kohnert2007constructing} \\ 
      627 & 307 & 16 & 396 & 204 & 188 & 627 & 396 & \cite{dissett2000combinatorial}, \cite{kohnert2007constructing} \\
      594 & 290 & 16 & 429 & 221 & 205 & 594 & 429 & \cite{dissett2000combinatorial}, \cite{kohnert2007constructing} \\
      561 & 273 & 16 & 462 & 238 & 222 & 561 & 462 & \cite{dissett2000combinatorial}, \cite{kohnert2007constructing} \\
      528 & 256 & 16 & 495 & 255 & 239 & 528 & 495 & \cite[Construction RT2]{calderbank1986geometry}, \cite{kohnert2007constructing} \\ 
      495 & 239 & 16 & 528 & 272 & 256 & 495 & 528 & complement \\
      462 & 222 & 16 & 561 & 289 & 273 & 462 & 561 & complement \\
      429 & 205 & 16 & 594 & 306 & 290 & 429 & 594 & complement \\
      396 & 188 & 16 & 627 & 323 & 307 & 396 & 627 & complement \\
      363 & 171 & 16 & 660 & 340 & 324 & 363 & 660 & complement \\
      330 & 154 & 16 & 693 & 357 & 341 & 330 & 693 & complement \\
      297 & 137 & 16 & 726 & 374 & 358 & 297 & 726 & complement \\
      264 & 120 & 16 & 759 & 391 & 375 & 264 & 759 & complement \\ 
      231 & 103 & 16 & 792 & 408 & 392 & 231 & 792 & complement \\
      198 &  86 & 16 & 825 & 425 & 409 & 198 & 825 & complement \\
      165 &  69 & 16 & 858 & 442 & 426 & 165 & 858 & complement does not exist \\
      132 &  52 & 16 & 891 & 459 & 443 & 132 & 891 & complement does not exist \\
       99 &  35 & 16 & 924 & 476 & 460 &  99 & 924 & complement does not exist \\
       66 &  18 & 16 & 957 & 493 & 477 &  66 & 957 & complement does not exist \\
       33 &   1 & 16 & 990 & 510 & 494 &  33 & 990 & complement does not exist \\
      \hline
    \end{tabular}  
    \caption{Feasible parameters for two-character sets in $\PG(9,2)$ excluding those for $\cP_i$, $\cS_i$, and their complements.}
    \label{table_feasible_parameters_sets_pg_9_2}
  \end{center}
\end{table}

\subsection{Feasible parameters of two-character sets in $\mathbf{\PG(10,2)}$}
In Table~\ref{table_feasible_parameters_sets_pg_10_2} we have listed the feasible parameters of two-character sets in $\PG(10,2)$ that are not covered by $\cS_i$ for 
$1\le i\le 10$, or their complements.

\begin{table}[htp]
  \begin{center}
    \begin{tabular}{|rrr|rrr|rr|l|}
      \hline
      $r$ & $\mu$ & $g$ & $n'$ & $s_0$ & $t_0$ & $a_s$ & $a_t$ & construction/non-existence \\
      \hline
      1780 & 884 & 16 &  712 & 360 & 328 & 1780 &  267 & \textbf{open} \\  
       276 & 132 & 16 &  759 & 407 & 375 &  276 & 1771 & \cite[Construction RT5]{calderbank1986geometry} \\
      1771 & 875 & 16 & 1288 & 648 & 616 & 1771 &  276 & complement \\
       267 & 123 & 16 & 1335 & 695 & 663 &  267 & 1780 & complement of open case \\
      \hline
      1335 & 663 & 32 &  267 & 139 & 123 & 1335 &  712 & \textbf{open} \\  
       759 & 375 & 32 &  276 & 148 & 132 &  759 & 1288 & \cite[Construction RT5]{calderbank1986geometry} \\ 
      1288 & 616 & 32 & 1771 & 891 & 875 & 1288 &  759 & complement \\
       712 & 328 & 32 & 1780 & 900 & 884 &  712 & 1335 & complement of open case \\
      \hline
    \end{tabular}  
    \caption{Feasible parameters for two-character sets in $\PG(10,2)$ excluding those for $\cS_i$ and their complements.}
    \label{table_feasible_parameters_sets_pg_10_2}
  \end{center}
\end{table}

We remark that a projective $32$-divisible binary code of length $712$ and a projective $16$-divisible binary code of length $267$ indeed exist, see \cite{kurz2021divisible}. However, 
the currently known constructions use more than two weights. 

\subsection{Feasible parameters of two-character sets in $\mathbf{\PG(11,2)}$}
In Table~\ref{table_feasible_parameters_sets_pg_11_2_1} and Table~\ref{table_feasible_parameters_sets_pg_11_2_2} we have listed the feasible parameters of two-character sets in $\PG(9,2)$ that are not covered by $\cP_i$ for $2\le i\le 63$, $\cS_i$ for 
$1\le i\le 11$, or their complements. There remains one block corresponding to the parameters of $\cK_i$ for $1\le i\le 62$. The parameters for $\cK_1,\dots,\cK_6$ as well as the complementary parameters for 
$\cK_{62},\dots,\cK_{57}$ cannot be attained:  
\begin{Lemma}
  \label{label_exl_k12_1}
  No projective $32$-divisible $[n,k]_2$-code exists for $n\in\{65,130,195,260,325,390\}$. 
\end{Lemma}
\begin{Proof}
  For $n\in \{65,130,195,260\}$ we apply \cite[Theorem 12]{honold2018partial}. For $n\in\{325,390\}$ we apply \cite[Corollary 6.15]{kurz2021divisible} with $t=5,6$.
\end{Proof}

\begin{table}[htp]
  \begin{center}
    \begin{tabular}{|rrr|rrr|rr|l|}
      \hline
      $r$ & $\mu$ & $g$ & $n'$ & $s_0$ & $t_0$ & $a_s$ & $a_t$ & construction/non-existence \\
      \hline 
      4030 & 2014 & 32 &   65 &   33 &    1 & 4030 &   65 & Lemma~\ref{label_exl_k12_1} \\   
      3965 & 1981 & 32 &  130 &   66 &   34 & 3965 &  130 & Lemma~\ref{label_exl_k12_1} \\
      3900 & 1948 & 32 &  195 &   99 &   67 & 3900 &  195 & Lemma~\ref{label_exl_k12_1} \\
      3835 & 1915 & 32 &  260 &  132 &  100 & 3835 &  260 & Lemma~\ref{label_exl_k12_1} \\ 
      3770 & 1882 & 32 &  325 &  165 &  133 & 3770 &  325 & Lemma~\ref{label_exl_k12_1} \\
      3705 & 1849 & 32 &  390 &  198 &  166 & 3705 &  390 & Lemma~\ref{label_exl_k12_1} \\
      3640 & 1816 & 32 &  455 &  231 &  199 & 3640 &  455 & \cite[Construction CY1]{calderbank1986geometry}, \cite{kohnert2007constructing} \\
      3575 & 1783 & 32 &  520 &  264 &  232 & 3575 &  520 & \textbf{open} \\    
      3510 & 1750 & 32 &  585 &  297 &  265 & 3510 &  585 & \textbf{open} \\
      3445 & 1717 & 32 &  650 &  330 &  298 & 3445 &  650 & \textbf{open} \\ 
      3380 & 1684 & 32 &  715 &  363 &  331 & 3380 &  715 & \textbf{open} \\
      3315 & 1651 & 32 &  780 &  396 &  364 & 3315 &  780 & \cite{bierbrauer1997family}, \cite{kohnert2007constructing} \\        -
      3250 & 1618 & 32 &  845 &  429 &  397 & 3250 &  845 & \cite{kohnert2007constructing} \\
      3120 & 1552 & 32 &  975 &  495 &  463 & 3120 &  975 & \cite{kohnert2007constructing} \\
      3055 & 1519 & 32 & 1040 &  528 &  496 & 3055 & 1040 & \cite{kohnert2007constructing} \\
      2990 & 1486 & 32 & 1105 &  561 &  529 & 2990 & 1105 & \cite{kohnert2007constructing} \\
      2925 & 1453 & 32 & 1170 &  594 &  562 & 2925 & 1170 & \cite{kohnert2007constructing} \\
      2860 & 1420 & 32 & 1235 &  627 &  595 & 2860 & 1235 & \textbf{open} \\     
      2795 & 1387 & 32 & 1300 &  660 &  628 & 2795 & 1300 & \cite{kohnert2007constructing} \\
      2730 & 1354 & 32 & 1365 &  693 &  661 & 2730 & 1365 & \cite[Construction CY1,CY2]{calderbank1986geometry}, \cite{kohnert2007constructing} \\
      2665 & 1321 & 32 & 1430 &  726 &  694 & 2665 & 1430 & \cite{kohnert2007constructing} \\
      2600 & 1288 & 32 & 1495 &  759 &  727 & 2600 & 1495 & \cite{kohnert2007constructing} \\ 
      2535 & 1255 & 32 & 1560 &  792 &  760 & 2535 & 1560 & \cite{kohnert2007constructing} \\
      2470 & 1222 & 32 & 1625 &  825 &  793 & 2470 & 1625 & \cite{kohnert2007constructing} \\
      2405 & 1189 & 32 & 1690 &  858 &  826 & 2405 & 1690 & \cite{kohnert2007constructing} \\
      2340 & 1156 & 32 & 1755 &  891 &  859 & 2340 & 1755 & \cite{kohnert2007constructing} \\
      2275 & 1123 & 32 & 1820 &  924 &  892 & 2275 & 1820 & \cite[Construction CY1,CY2]{calderbank1986geometry}, \cite{kohnert2007constructing} \\ 
      2210 & 1090 & 32 & 1885 &  957 &  925 & 2210 & 1885 & \cite{kohnert2007constructing} \\
      2145 & 1057 & 32 & 1950 &  990 &  958 & 2145 & 1950 & \cite{kohnert2007constructing} \\
      2080 & 1024 & 32 & 2015 & 1023 &  991 & 2080 & 2015 & \cite[Construction RT2]{calderbank1986geometry}, \cite{kohnert2007constructing} \\
      2015 &  991 & 32 & 2080 & 1056 & 1024 & 2015 & 2080 & complement \\
      1950 &  958 & 32 & 2145 & 1089 & 1057 & 1950 & 2145 & complement \\
      1885 &  925 & 32 & 2210 & 1122 & 1090 & 1885 & 2210 & complement \\
      1820 &  892 & 32 & 2275 & 1155 & 1123 & 1820 & 2275 & complement \\
      1755 &  859 & 32 & 2340 & 1188 & 1156 & 1755 & 2340 & complement \\
      1690 &  826 & 32 & 2405 & 1221 & 1189 & 1690 & 2405 & complement \\
      1625 &  793 & 32 & 2470 & 1254 & 1222 & 1625 & 2470 & complement \\
      1560 &  760 & 32 & 2535 & 1287 & 1255 & 1560 & 2535 & complement \\
      1495 &  727 & 32 & 2600 & 1320 & 1288 & 1495 & 2600 & complement \\
      1430 &  694 & 32 & 2665 & 1353 & 1321 & 1430 & 2665 & complement \\
      1365 &  661 & 32 & 2730 & 1386 & 1354 & 1365 & 2730 & complement \\
      1300 &  628 & 32 & 2795 & 1419 & 1387 & 1300 & 2795 & complement \\
      1235 &  595 & 32 & 2860 & 1452 & 1420 & 1235 & 2860 & complement of open case \\
      1170 &  562 & 32 & 2925 & 1485 & 1453 & 1170 & 2925 & complement \\
      1105 &  529 & 32 & 2990 & 1518 & 1486 & 1105 & 2990 & complement \\
      1040 &  496 & 32 & 3055 & 1551 & 1519 & 1040 & 3055 & complement \\
       975 &  463 & 32 & 3120 & 1584 & 1552 &  975 & 3120 & complement \\
       845 &  397 & 32 & 3250 & 1650 & 1618 &  845 & 3250 & complement \\
       780 &  364 & 32 & 3315 & 1683 & 1651 &  780 & 3315 & complement \\
      \hline
    \end{tabular}  
    \caption{Feasible parameters for two-character sets in $\PG(11,2)$ excluding those for $\cP_i$, $\cS_i$, and their complements -- part 1.}
    \label{table_feasible_parameters_sets_pg_11_2_1}
  \end{center}
\end{table}

\begin{table}[htp]
  \begin{center}
    \begin{tabular}{|rrr|rrr|rr|l|}
      \hline
      $r$ & $\mu$ & $g$ & $n'$ & $s_0$ & $t_0$ & $a_s$ & $a_t$ & construction/non-existence \\
      \hline 
       715 &  331 & 32 & 3380 & 1716 & 1684 &  715 & 3380 & complement of open case \\
       650 &  298 & 32 & 3445 & 1749 & 1717 &  650 & 3445 & complement of open case \\
       585 &  265 & 32 & 3510 & 1782 & 1750 &  585 & 3510 & complement of open case \\
       520 &  232 & 32 & 3575 & 1815 & 1783 &  520 & 3575 & complement of open case \\
       455 &  199 & 32 & 3640 & 1848 & 1816 &  455 & 3640 & complement \\
       390 &  166 & 32 & 3705 & 1881 & 1849 &  390 & 3705 & complement does not exist \\
       325 &  133 & 32 & 3770 & 1914 & 1882 &  325 & 3770 & complement does not exist \\
       260 &  100 & 32 & 3835 & 1947 & 1915 &  260 & 3835 & complement does not exist \\
       195 &   67 & 32 & 3900 & 1980 & 1948 &  195 & 3900 & complement does not exist \\                                                                                                                                
       130 &   34 & 32 & 3965 & 2013 & 1981 &  130 & 3965 & complement does not exist \\
        65 &    1 & 32 & 4030 & 2046 & 2014 &   65 & 4030 & complement does not exist \\
      \hline
      2808 & 1400 & 64 &  234 &  122 &  106 & 2808 & 1287 & \cite[Construction FE3]{calderbank1986geometry}, \cite{chen2006constructions}, \cite{kohnert2007constructing} \\      
      1400 &  696 & 64 &  245 &  133 &  117 & 1400 & 2695 & \textbf{open} \\       
      2295 & 1143 & 64 &  270 &  142 &  126 & 2295 & 1800 & \cite{bierbrauer1997family}, \cite{kohnert2007constructing} \\ 
      1911 &  951 & 64 &  273 &  145 &  129 & 1911 & 2184 & \cite{chen2006constructions}, \cite{kohnert2007constructing} \\        
      2184 & 1032 & 64 & 3822 & 1918 & 1902 & 2184 & 1911 & complement \\
      1800 &  840 & 64 & 3825 & 1921 & 1905 & 1800 & 2295 & complement \\ 
      2695 & 1287 & 64 & 3850 & 1930 & 1914 & 2695 & 1400 & complement of open case \\
      1287 &  583 & 64 & 3861 & 1941 & 1925 & 1287 & 2808 & complement \\ 
      \hline
      3185 & 1585 & 32 &  910 &  462 &  430 & 3185 &  910 & \cite[Construction CY1,CY2]{calderbank1986geometry} \\
       910 &  430 & 32 & 3185 & 1617 & 1585 &  910 & 3185 & complement \\ 
      \hline
      3861 & 1925 & 16 & 1287 &  647 &  583 & 3861 &  234 & \cite{kohnert2007constructing}\footnote{At http://www.tec.hkr.se/$\sim$chen/research/2-weight-codes, refering to \cite{kohnert2007constructing}, an explicit generator matrix is given. However, in \cite{kohnert2007constructing} those parameters are not treated and the corresponding webpage http://linearcodes.uni-bayreuth.de/twoweight is not accesible any more. Each of the examples for $n'=234$ can also be used to construct such a code.} \\
       245 &  117 & 16 & 1400 &  760 &  696 &  245 & 3850 & \textbf{open} \\ 
      3825 & 1905 & 16 & 1800 &  904 &  840 & 3825 &  270 & \cite{bierbrauer1997family}, \cite{kohnert2007constructing} \\  
       273 &  129 & 16 & 1911 & 1015 &  951 &  273 & 3822 & \cite{de1995some}, \cite{kohnert2007constructing} \\
      3822 & 1902 & 16 & 2184 & 1096 & 1032 & 3822 &  273 & complement \\ 
       270 &  126 & 16 & 2295 & 1207 & 1143 &  270 & 3825 & complement \\ 
      3850 & 1914 & 16 & 2695 & 1351 & 1287 & 3850 &  245 & complement of open case \\
       234 &  106 & 16 & 2808 & 1464 & 1400 &  234 & 3861 & complement \\ 
      \hline
    \end{tabular}  
    \caption{Feasible parameters for two-character sets in $\PG(11,2)$ excluding those for $\cP_i$, $\cS_i$, and their complements -- part 2.}
    \label{table_feasible_parameters_sets_pg_11_2_2}
  \end{center}
\end{table}
W.r.t.\ the series $\cK_i$ we remark that it is currently unknown whether projective $32$-divisible $[n,k]_2$-codes exist for $n\in\{520,585,650,715\}$, see e.g.\ \cite{kurz2021divisible}. 
W.r.t.\ the other open cases no construction for a projective $\Delta$-divisible $[n,k]_2$-code is known for $(\Delta,n)=(16,245)$ and $(\Delta,n)=(64,1400)$.

\end{document}